\documentclass[a4paper,12pt]{amsart}

\usepackage{setspace}

\usepackage{amsmath,amssymb}
\usepackage{amsfonts}
\usepackage{epsfig}

\def\n{\noindent}
\def\cB{{\mathcal B}}
\def\cD{{\mathcal D}}
\def\cE{{\mathcal E}}
\def\cH{{\mathcal H}}
\def\cL{{\mathcal L}}
\def\cN{{\mathcal N}}
\def\cX{{\mathcal X}}

\def\ZZ{{\mathbb Z}}

\def\RR{{\mathbb R}}
\def\PP{{\mathbb P}}
\def\Pr#1{\PP\left\{#1\right\}}
\def\EE{{\mathbb E}}
\def\Ex#1{{\EE\left\{#1\right\}}}

\def\Ltwo{L^2(\RR )}

\def\tsigma{{\tilde\sigma}}
\def\ttheta{{\tilde\theta}}
\def\tpi{{\tilde\pi}}
\def\tP{{\tilde P}}
\def\tp{{\tilde p}}
\def\tepsilon{{\tilde \epsilon}}

\def\dd{\displaystyle}

\newtheorem{lemma}{Lemma}
\newtheorem{proposition}{Proposition}
\newtheorem{definition}{Definition}

\newtheorem{corollary}{Corollary}
\newtheorem{remark}{Remark}

\begin{document}
\doublespacing
%\onehalfspacing

\title{A Hybrid Scheme for Encoding Audio Signal using Hidden Markov Models
of Waveforms}

\author{S. MOLLA and B. TORRESANI}

\address{Laboratoire d'Analyse, Topologie et Probabilit\'es,
CMI, Universit\'e de Provence,
39 rue F. Joliot-Curie, 13453 Marseille Cedex 13,
France.
}
\email{molla@cmi.univ-mrs.fr ; torresan@cmi.univ-mrs.fr}

\begin{abstract}
This paper reports on recent results related to
audiophonic signals encoding using time-scale and
time-frequency transform. More precisely,
non-linear, structured approximations for tonal
and transient components using local cosine and wavelet
bases will be described, yielding expansions of audio
signals in the form
{\em tonal + transient + residual}.
We describe a general formulation involving hidden Markov models,
together with corresponding rate estimates.
Estimators for the balance transient/tonal are also discussed.
\end{abstract}

\maketitle
\markboth{\qquad\em S. Molla and B. Torr\'esani
%\qquad\qquad\qquad\qquad
\hfill submitted to ACHA, revised \today}
{{\em An Hybrid Audio Scheme using Hidden Markov Models of Waveforms.\hfill submitted to ACHA},\qquad}

\section{Introduction: structured hybrid models}
\label{se:intro}
Recent signal processing studies have shown the importance of
sparse representations for various tasks, including signal
and image compression (obviously), de-noising, signal
identification/detection,... Such sparse
representations are generally achieved using suitable orthonormal bases
of the considered signal space. However, recent developments
also indicate that redundant systems, such as frames, or more
general ``waveform dictionaries'' may yield substantial
gains in this context, provided that they are
sufficiently adapted to the signal/image to be described.

From a different point of view, it has also been shown by several
authors that in a signal or image compression context,
significant improvements may be achieved by introducing
{\em structured approximation} schemes, namely schemes
in which structured sets of coefficients are considered rather than
isolated ones.

The goal of this paper is to describe a new approach that implements
both ideas, via a hybrid model involving sparse, structured,
random wavelet/MDCT\footnote{MDCT: Modified Discrete Cosine Transform.}
expansions, where the sets of considered
coefficients (the {\em significance maps}) are described
via suitable (hidden) Markov models.

This work is mainly motivated by audio coding applications, to which
we come back after describing the models and corresponding estimation
algorithms. However, similar ideas may clearly be developed in different
contexts, including image~\cite{Meyer02multilayered}
and image sequence coding, where both ingredients
(hybrid and structured models) have already been exploited.
\subsection{Generalities, sparse expansions in redundant systems}
Very often, signals turn out to be made of several
components, of significantly different nature. This is the
case for ``natural images'', which may contain
edge information, regular textures, and ``non-stationary''
textures (which carry 3D information.) This is also the case for
audio signals, which among other features, contain transient and
tonal components~\cite{Daudet02hybrid}, on which we shall focus
more deeply. It is known that such different features may be
represented efficiently in specific orthonormal bases. Following the
philosophy of transform coding, this suggests
to consider redundant systems made out by concatenation of
several families of bases. Such systems have been considered for example
in~\cite{Donoho01uncertainty,Elad02generalized,Gribonval03sparse},
where the problem of selecting the ``sparsest'' expansion
through linear programing has been considered.

Focusing on the particular application to audio signals, and
limiting ourselves to transient and tonal features, we are
naturally led to consider a generic redundant dictionary
made out of two orthonormal bases, denoted by $\psi_\lambda$
and $w_\delta$ respectively (typically a wavelet and an MDCT basis),
and signal expansions of the form
\begin{equation}
\label{fo:sig.model}
x = \sum_{\lambda\in\Lambda} \alpha_{\lambda}\psi_{\lambda} +
\sum_{\delta\in\Delta} \beta_{\delta} w_{\delta} +r\ ,
\end{equation}
where $\Lambda$ and $\Delta$ are (small, and this will be the main
{\em sparsity} assumption) subsets of the index sets,
hereafter termed {\em significance maps}.
The nonzero coefficients $\alpha_\lambda$
are independent $\cN(0,\sigma_\lambda^2)$
random variables, and the nonzero coefficients $\beta_\delta$
are independent $\cN(0,\tsigma_\delta^2)$ random variables:
$r$ is a residual signal, which is not sparse with respect
to the two considered bases (we shall talk of {\em spread residual}),
and is to be neglected or described differently.

The approach developed
in~\cite{Donoho01uncertainty,Elad02generalized,Gribonval03sparse}
may be criticized in several respects when it comes to
practical implementation in a coding perspective. On one hand,
it is not clear that the corresponding linear programing
algorithms are compatible with practical constraints, in terms of
CPU and memory requirements\footnote{for example, for audio signals typically
sampled at 44.1 kHz.}. Also, models exploiting solely sparsity
arguments cannot capture one of the main features of some
signal classes, namely the {\em persistence} property:
significant coefficients have a tendency to form ``clusters'',
or ``structured sets''. For example, in
an audio coding context, the significance maps
take the form of {\em ridges} (i.e. ``time-persistent'' sets,
see e.g.~\cite{Carmona98practical,Delprat92asymptotic} in a different context)
for the MDCT map $\Delta$, and binary trees for the wavelet map $\Lambda$.
This remark has been exploited in various instances,
for example in the context of the sinusoidal models for
speech~\cite{Macaulay86speech}, of for image
coding~\cite{Cohen01tree,Crouse98wavelet,Said96new,Shapiro93embedded}

Several models may be considered for the $\Lambda$ and $\Delta$
sets (termed {\em significance maps}), with variable levels
of complexity. If only sparsity is used, they may be chosen
uniformly distributed (in a finite dimensional context.)
We shall rather work in a more complex context, and use
(hidden) Markov chains to describe the MDCT ridges in $\Delta$
(in the spirit of the sinusoidal models of speech),
and (hidden) binary Markov trees for the wavelet map $\Lambda$,
following~\cite{Crouse98wavelet}. This not only yields a better
modeling of the features of the signal, but also provides
corresponding estimation algorithms.

To be more specific, a tonal signal is modeled as
$$
x_{ton} = \sum_{\delta\in\Delta} \beta_\delta w_\delta\ ,
$$
%$x_{nton}$ being some non-tonal residual signal, and
the functions $w_\delta$ being local cosine functions.
The (significant) coefficients $\beta_\delta$, $\delta\in\Delta$
are $\cN(0,\tsigma_\delta^2)$ independent random variables.
The index $\delta$ is in fact a pair of time-frequency indices
$\delta=(k,\nu)$, and the significance map $\Delta$ is
characterized by a ``fixed frequency'' Markov chain
(see e.g.~\cite{Karlin97first} for a simple account), hence by
a set of initial frequencies $\nu_1,\dots\nu_N$ and transitions matrices
$\tP_1,\dots\tP_N$ (one for each frequency bin).

Globally, the tonal model is characterized by the set of matrices
$\tP_n$, and the variances $\tsigma_\delta^2$ of the two states,
which are assumed to be time invariant, and on which additional
constraints may be imposed. The tonal model is described in some
details in section~\ref{se:tonal}.

A similar model, using Hidden Markov {\em trees} of wavelet
coefficients~\cite{Crouse98wavelet} may be develop to
describe the transient layer in the signal:
$$
x_{tr} = \sum_{\lambda\in\Lambda} \alpha_\lambda \psi_\lambda\ ,
$$
$\psi$ being a wavelet with good time localization.
The rationale is now to
model the {\em scale persistence} of large wavelet coefficients
of the transients, exploiting the intrinsic dyadic tree
structure of wavelet coefficients (see {\sc Figure}~\ref{fi:dyadtree} below.)
Again, the significant wavelet coefficients
$\{\alpha_\lambda,\lambda\in\Lambda\}$ of the signal
are modeled as independent $\cN(0,\sigma_\lambda^2)$ random
variables.
The index $\lambda$ is in fact a pair of scale-time indices
$\delta=(j,k)$, and the significance map $\Lambda$ is
characterized by a ``fixed time'' Markov chain, hence by
corresponding ``scale to scale'' transition matrices $P_j$
(with additional constraints which ensure that significant
coefficients inherit a natural tree structure, see below.)
%, i.e. numbers of the type
%$$
%\pi_j = \Pr{(j+1,2k)\in\Lambda \vert (j,k)\in\Lambda}
%= \Pr{(j+1,2k+1)\in\Lambda \vert (j,k)\in\Lambda}\ .
%$$
%For practical reasons, it is convenient to impose the
%further constraint $\pi_j'=1$, $\pi_j'$ being the
%scale-persistence probability of ``residual'' states.
%This has the advantage of ensuring the connectedness
%of the binary trees, which is useful for significance map
%encoding, and probably fairly realistic in terms of
%transient modeling.

The transient model is therefore characterized by the variances
of wavelet coefficients in $\Lambda$ and $\Lambda^c$, and the
persistence probabilities, for which estimators may be constructed.
The transient states estimation itself is also performed
via classical methods.
These aspects are described in section~\ref{se:transient}.

\subsection{Recursive estimation}
Several approaches are possible to estimate the significance maps
and corresponding coefficients in models such as~(\ref{fo:sig.model}),
ranging from the above mentioned linear programing schemes
(see for example~\cite{Chen98atomic}) to greedy algorithms,
including for instance Matching
pursuit~\cite{Gribonval99approximation,Mallat93matching}. The
procedure we use is in some sense intermediate between these two
extremes, in the spirit of the techniques used in~\cite{Berger94removing}.
We consider a dictionary made of two (orthonormal) bases; 
a first layer is estimated, using the first basis, and a second
layer is estimated from the residual, using the second basis.
The main difficulty of such an approach lies on the fact that
the number of significant elements from the first basis has to
be known in advance (or at least estimated.) In other terms,
the cardinalities $|\Lambda|$ and $|\Delta|$ of the significance
maps have to be known. This is important, since an underestimation
or overestimation of $|\Delta|$ (assuming that the $\Delta$-layer
is estimated first) will ``propagate'' to the estimation of the
second layer (the $\Lambda$-layer.)

In the framework of the the Gaussian random sparse models
studied below, it is possible to derive a priori estimates
for the cardinalities $|\Lambda|$ and $|\Delta|$, using
information measures in the spirit of those proposed
in~\cite{Wickerhauser94adapted} and studied
in~\cite{Trgo95relation}. Consider the geometric means
of estimated $\psi_\lambda$ and $w_\delta$ coefficients
\begin{equation}
\label{fo:N}
\hat N_\psi = \left(\prod_{n=1}^N |\langle x,\psi_n\rangle|^2\right)^{1/N}\ 
\hbox{and }
\qquad
\hat N_w = \left(\prod_{n=1}^N |\langle x,w_n\rangle|^2\right)^{1/N}\ .
\end{equation}
Then, assuming spartity, the indices
\begin{equation}
I_{w} = \frac{\hat N_{\psi}}{\hat N_{\psi}+ \hat N_{w}}\ ;\qquad
I_{\psi} = \frac{\hat N_{w}}{\hat N_{\psi}+\hat  N_{w}}\ ,
\end{equation}
turn out to provide estimates for the proportion of significant
$w$ and $\psi$ coefficients. The rationale is the fact that under
sparsity assumptions (i.e. if $\Delta$ and $\Lambda$ are small enough),
most coefficients $\langle x,\psi_n\rangle$ (resp.
$\langle x,w_n\rangle$) will come from the tonal (resp. transient)
layer of the signal, and therefore give information about it.
This aspect is discussed in more details in section~{\ref{se:balance}.

\subsection{Audio coding applications}
As mentioned earlier, the primary motivation for this work was
audio coding. We briefly sketch here the assets of the model we are
developing in such a context.

Coding involve (lossy) quantization of the selected coefficients
$\{\langle x,w_\delta\rangle,\delta\in\Delta\}$ and
$\{\langle x,\psi_\lambda\rangle,\lambda\in\Lambda\}$. These
are Gaussian random variables, which means that corresponding
rate and distortion estimates may be obtained. We notice that the
introduction of {\em structured} significance maps does {\em not} improve
the quality of the approximation (as measured by $L^2$ distortion);
however it improves the efficiency of significance map encoding
(see below); in addition, for audio applications, since structured
significance maps seem to be relevant in the sense that they describe
more accurately elementary ``sound objects'', they often yield
better approximations of audio signals, from perceptual points of
view.

Besides coefficients, significance maps have to be encoded as
well. However, the Markov models make it possible to compute explicitly
the probabilities of ridges lengths (for $\Delta$)
and trees sizes, which allows one to obtain directly
the corresponding optimal lossless code. Again,
rate estimates may be derived explicitly.

It is also worth pointing out some important issues
(in a coding perspective), which we shall not address here.
The first one is the encoding of the residual signal
$$
x_{res} = x - x_{ton} - x_{tr}\ .
$$
It was suggested in~\cite{Daudet02hybrid} that the residual
may be encoded using standard LPC techniques\footnote{LPC: Linear
Predictive Coding.}.
However, it appears that in most situations (at least for large enough
bit rates), encoding the residual is not necessary,
the transient and tonal layers providing a satisfactory
description of the signal.

A second point is related to the implementation of
perceptive arguments (e.g. masking): the goal is not really to
obtain a lossy description of the signal with a {\em small}
distortion: the distortion is rather expected to be inaudible,
which has little to do with its $\ell^2$ norm. In the
proposed scheme, this aspect will be addressed at the level of
coefficient quantization (as in most perceptive coders.) However
let us point out that the ``structural decomposition''
involving well defined tonal and transient layers shall make it
possible to implement separately frequency masking on the
tonal layer, and time masking on the transient layer, which is
a completely original approach. This work (in progress) will be
partly reported in~\cite{Daudet03HSAC}.

%
%%%%
%
\section{Structured Markov model for tonal}
\label{se:tonal}
We start with a description of the first layer of the model.
We make use of the local cosine bases constructed by Coifman
and Meyer. Let us briefly recall here the construction, in the
case we shall be interested in here. Let $\ell\in\RR^+$ and
$\eta\in\RR^+,\eta<\ell/2$. Let $w$ be a smooth function
(called the basic window) satisfying the following properties:
\begin{eqnarray}
\hbox{supp}(w)&\subset& [0-\eta,\ell+\eta]\\
w(-\tau)&=& w(\tau)\quad\hbox{for all }|\tau|\le\eta\\
w(\ell-\tau)&=& w(\ell+\tau)\quad\hbox{for all }|\tau|\le\eta\\
\sum_{k} w(t-k\ell)^2 &=&1\ ,\quad\forall t\ .
\end{eqnarray}
and set
\begin{equation}
w_{kn}(t) = \sqrt{\frac{2}\ell}\ w(t-k\ell)\ 
\cos\left(\frac{\pi(n+1/2)}\ell (t-k\ell)\right)\ ,\quad
n\in\ZZ^+,\ k\in\ZZ\ .
\end{equation}
Then it may be proved that the collection of such functions,
when $n$ spans $\ZZ^+$ and $k$ spans $\ZZ$, forms an
orthonormal basis of $\Ltwo$. Versions adapted to spaces
of functions with bounded support, as well as discrete versions, may
also be obtained easily. We refer to~\cite{Wickerhauser94adapted}
for a detailed account of such constructions.
The classical choice for such functions amounts to take an arc
of sine wave as function $w$. We shall limit ourselves to the
so-called ``maximally smooth'' windows, by setting $\eta =\ell/2$.

\begin{figure}
\centerline{
\includegraphics[height=7cm,width=13cm]{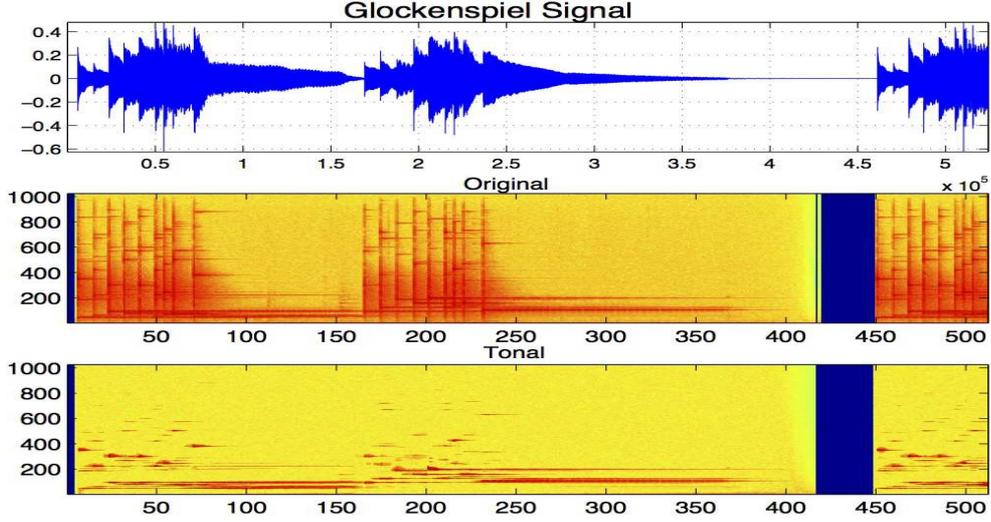}
}
\caption{Estimating a tonal layer; top: glockenspiel signal;
middle: logarithm of absolute value of MDCT coefficients of the signal;
bottom: logarithm of absolute value of MDCT coefficients of a tonal
layer, estimated using ``horizontal'' structures in MDCT coefficients.}
\label{fi:gspi.full}
\end{figure}

In the framework of the recursive estimation scheme we are about to
describe, the simplest (and natural) idea would be to start by
expanding the signal with respect to a local cosine basis, and pick the
largest coefficients (in absolute value, after appropriate weighting
if needed) to form a best $N$-term approximation~\cite{Daudet02hybrid}.
However, as may be seen in the middle image of {\sc Figure}~\ref{fi:gspi.full},
such a strategy would automatically ``capture'' local cosine coefficients
which definitely belong to transients (i.e. seem to form localized,
``vertical'' structures.)
In order to avoid capturing such undesired coefficients,
it is also natural to use the ``structure'' of MDCT coefficients of
tonals, i.e. the fact that they have a tendancy to form ``horizontal ridges''.
This is the purpose of the tonal model described below. In the glockenspiel
example of {\sc Figure}~\ref{fi:gspi.full}, such a strategy produces
a tonal layer whose MDCT is exhibited in the bottom image, from which it
is easily seen that only ``horizontally structured'' coefficients have been
retained.

\subsection{Model and consequences}
\label{sub:ton.model}
In the framework of the recursive approach, the signal
is modeled as a {\em structured harmonic mixture of Gaussians},
i.e. expanded into an MDCT basis, with given cutoff frequency
$N$
\begin{equation}
x = \sum_{n=0}^{N-1}\sum_{k} Y_{kn} w_{kn}\ ,
\end{equation}
where the coefficients of the expansion are (real,
continuous) random variables $Y_{kn}$ whose distribution
is governed by a family of ``fixed frequency''
Hidden Markov chains (HMC) $X_{kn}, k=1,\dots$.
According to the usual practice, we shall denote by
$Y_{1:k,n}$ (resp. $X_{1:k,n}$) the random vector
$(Y_{1n},\dots Y_{kn})$ (resp. $(X_{1n},\dots X_{kn})$),
and use a similar notation for the corresponding values
$(y_{1n},\dots y_{kn})$ (resp. $(x_{1n},\dots x_{kn})$.)
$\rho_{Y_{1:k,n}}$ and $\rho_{Y_k}$ will denote the
joint density of $Y_{1:k,n}$ and the density of $Y_{kn}$
respectively, and the density of $Y_{kn}$ conditioned by
$X_{kn}$, assumed to be independent of $k$,
will be denoted by
$$
\psi_{n}(y\vert x) = \rho_{Y_{kn}}(y\vert X_{kn}=x)\ ,\quad x=T,N\ .
$$
To be more precise, the model is characterized as follows:
\begin{enumerate}
\item[{\em i.} ] For all $n$, $X_{\cdot n}$ is a Markov chain with state space
$$
\cX = \left\{T,R\right\}
$$
(``tonal'' and ``residual'', or non-tonal)
and transition matrix $\tP_n$, of the form
$$
\tP_n =\begin{pmatrix}
\tpi_n&1 -\tpi_n\\
1-\tpi_n'&\tpi_n'
\end{pmatrix}
$$
the numbers $\tpi_n,\tpi_n'$ being the {\em persistence probabilities}
of the tonal and residual states: for all $n$
%$$
%\tpi = \Pr{(k,n)\in\Delta \vert (k-1,n)\in\Delta}\ ,\quad
%\tpi' = \Pr{(k,n)\not\in\Delta \vert (k-1,n)\not\in\Delta}\ .
%$$
\begin{eqnarray}
\tpi_n = \Pr{X_{kn}=T \vert X_{k-1\,n}=T}\ ,\\
\tpi_n' = \Pr{X_{kn}=R \vert X_{k-1\,n}=R}\ .
\end{eqnarray}
The initial frequencies of $T$ and $R$ states will be denoted by
$\nu_n$ and $1-\nu_n$ respectively. For the sake of simplicity,
we shall generally assume that the initial frequencies coincide
with the equilibrium frequencies of the chain:
$$
\nu_n^{(e)} = \frac{1-\tpi_n'}{2-\tpi_n-\tpi_n'}\ ,
$$
%\item For all $n$, the initial frequencies of the Markov chain $X^{(n)}$
%are the equilibrium frequencies of $\tP_n$.
\item[{\em ii.} ] The (emitted) coefficients $Y_{kn}$ are continuous random
variables, with densities denoted by $\rho_{Y_{1:k\,n}}(y_{1:k\,n})$,
\item[{\em iii.} ] The distribution of the (emitted)
coefficients $Y_{kn}$ depends
only on the corresponding hidden state $X_{kn}$; for each
$n$, the coefficients $Y_{kn}$ are independent conditional
to the hidden states, and their distribution do not depend
on the time index $k$ (but does depend on the frequency index $n$.)
We therefore denote
$$
\rho_{Y_{1:k n}}\left(y_{1:k n}\vert X_{1:k n}=x_{1:k n}\right)
= \prod_{i=1}^k \psi_{n}(y_{in}\vert x_{in})\ ,
$$
\item[{\em iv.} ] In order to model audio signal, we shall
limit ourselves to centered gaussian models for the
densities $\psi_k$. The latter 
are therefore completely determined by their variances: a large
variance $\sigma_T^2$ for the $T$ type coefficients, and a small
variance $\sigma_R^2$ for the $R$ type coefficients.
\end{enumerate}
Therefore, the model is completely characterized by the parameter set
$$
\ttheta = \{\tpi_n,\tpi_n';\nu_n; \tsigma_{T,n},\tsigma_{R,n};\,
n=0,\dots N-1\}\ .
$$
Given these parameters, one may compute explicitly the likelihood of
any configuration of coefficients. Using ``routine'' HMC
techniques, it is also possible to obtain explicit formulas for the
likelihood of any hidden states configuration, conditional to the
coefficients. We refer to~\cite{Rabiner89tutorial} for a detailed
account of these aspects.

\begin{remark}\rm
Notice that in this version of the model, the transition matrix $\tP$
is assumed to be frequency independent. More general models involving
frequency dependent $\tP$ matrices (or further generalizations) may
be constructed, without much modifications of the overall approach.
\end{remark}
Given a signal model as above, we may define the tonal layer of such a signal.
\begin{definition}
\label{def:ton.layer}
Let $x$ be signal modeled as a hidden Markov chain MDCT as above,
and let 
\begin{equation}
\label{def:ton.map}
\Delta = \left\{(k,n)\vert X_{kn}=T\right\}\ .
\end{equation}
$\Delta$ is called the {\em tonal significance map}
of $x$. Then the tonal and non tonal layers
are given by
\begin{eqnarray}
\label{fo:ton.part}
x_{ton} &=& \sum_{\delta\in\Delta} \beta_{\delta}
w_{\delta}\ ,\\
\label{fo:nton.part}
x_{nton} &=& x - x_{ton}
\end{eqnarray}
\end{definition}
This definition makes it possible to obtain simple estimates
for quantities of interest, such as the energy of a tonal
signal, or the number of MDCT coefficients needed to encode it.
For example, considering a time frame of $K$ consecutive
windows (starting from $k=0$ for simplicity\footnote{In fact,
this choice of origin matters only if the initial frequency
$\nu$ of the chain is not assumed to equal the equilibrium
frequency $\nu^{(e)}$, which will not be the case in the
situations we consider.}), and a frequency domain
$\{0,\dots N-1\}$, we set
$$
\Delta_{(K,N)} = \Delta\cap \left(\{0,\dots K-1\}
\times\{0,\dots N-1\}\right)\ ,
$$
and we denote by
\begin{eqnarray}
\tilde N_n^{(K)} &=& |\Delta_{(K,\{n\})}|\\
\tau_n^{(K)} &=& \Ex{\frac{\tilde N_n^{(K)}}K}
\end{eqnarray}
the random variables describing respectively the number and the
expected proportion of $T$ type coefficients in the frequency bin $n$,
within a time frame of $K$ consecutive windows.
\begin{proposition}
\label{prop:ton.rate}
With the notations of Definition~\ref{def:ton.layer}, 
the average proportion of $T$ type coefficients within the
time frame $\{0,\dots K-1\}$ in the frequency
bin $n$ is given by
\begin{equation}
\tau_n^{(K)} =
%\Ex{\frac{|\tilde N^{(K)}_n}K}\\
\frac1{K(2-\tpi_n-\tpi_n')}
\left[\nu_n ((\!\tpi_n\!+\!\tpi_n'\!-\!1\!)^K\!)\! +\! (\!1\!-\!\tpi_n'\!)
\left(K - \frac{1\!-\!(\!\tpi_n\!+\!\tpi_n'\!-\!1\!)^K}
{2\!-\!\tpi_n\!-\!\tpi_n'}\right)
\right]
\end{equation}
\end{proposition}
\noindent{\em Proof:} From classical properties of HMC, we have that
$$
\begin{pmatrix}
\Pr{X_{k\,n}=T}\\
\Pr{X_{k\,n}=R}
\end{pmatrix}
= \left(\tP^t\right)^k
\begin{pmatrix}
\nu_n\\1-\nu_n
\end{pmatrix}\ ,
$$
the superscript ``$t$'' denoting matrix transposition.
After some algebra, we obtain the following expressions:
\begin{eqnarray}
\nonumber
\Pr{X_{k\,n}=T} &=& \frac{
\left((1-\tpi_n)\nu_n - (1-\tpi_n')(1-\nu_n)\right)
\left(\tpi_n+\tpi_n'-1\right)^k
+ (1-\tpi_n')}
{2-\tpi_n-\tpi_n'}\\
\nonumber
&=& \nu_n \left(\tpi_n+\tpi_n'-1\right)^k + 
\frac{1-\tpi_n'}{2-\tpi_n-\tpi_n'}\,
\left(1-(\tpi_n+\tpi_n'-1)^k\right)\ .
\end{eqnarray}
Similarly, we obtain for $\Pr{X_{k\,n}=R}=1-\Pr{X_{k\,n}=T}$
$$
\Pr{X_{k\,n}=R} = (1-\nu_n) \left(\tpi_n+\tpi_n'-1\right)^k
+ \frac{1-\tpi_n}{2-\tpi_n-\tpi_n'}\,\left(1-(\tpi_n+\tpi_n'-1)^k\right)\ .
$$
Finally, the result is obtained by replacing $ \Pr{X_{k\,n}=T}$
with its expression in
$$
\Ex{\tilde N^{(K)}_n} = \sum_{k=0}^{K-1} \Pr{X_{k\,n}=T}\ ,
$$
which yields the desired expression.
\qed

Notice that in the limit of large time frames, one obtains the simpler
estimate
$$
\lim_{K\to\infty}\tau_n^{(K)}
%\Ex{\frac{|\Delta\cap\{k,\dots k+K-1\}||}K} 
=\frac{1-\tpi_n'}{2-\tpi_n-\tpi_n'} = \nu_n^{(e)}\ ,
$$
which of course does not depend any more on $K$.

The energy of the tonal layer is also completely characterized
by the parameters of the model, and has a simple behavior.
\begin{proposition}
\label{prop:ton.energy}
With the same notations as before, conditional to the parameters of the
model, we have
\begin{equation}
\begin{split}
\Ex{\frac1{K}\sum_{\delta\in\Delta_{(K,N)}} |Y_\delta|^2} =&
\frac1{K}\sum_{n=0}^{N-1}\frac1{2-\tpi_n-\tpi_n'}\bigg[
\left(1-(\tpi_n+\tpi_n'-1)^K\right)\nu_n\tsigma_{T,n}^2\\
&\quad+ (1-\tpi_n')\left(K-\frac{1-(\tpi_n+\tpi_n'-1)^K}{2-\tpi_n-\tpi_n'}
\right)\tsigma_{T,n}^2\bigg]
\end{split}
\end{equation}
\end{proposition}
\noindent{\em Proof:} the result follows from the fact that
conditional to the hidden states, the considered random variables
at fixed frequency are i.i.d. $\cN(0,\sigma_{T,n}^2)$ random variables.
It is then enough to plug the expression of $\tau_n^{(K)}$ obtained
above in the $L^2$ norm of the tonal layer.\qed

Again, the latter expression simplifies in the limit $K\to\infty$,
or if the initial frequencies of the
chains $X_n$ are assumed to equal the equilibrium frequencies. In that
situation, we obtain
\begin{equation}
\lim_{K\to\infty}\Ex{\frac1{K}\sum_{\delta\in\Delta_{(K,N)}} |Y_\delta|^2} = 
\sum_{n=0}^{N-1} \frac{1-\tpi_n'}{2-\tpi_n-\tpi_n'}\,\tsigma_{T,n}^2
= \sum_{n=0}^{N-1} \nu_n^{(e)}\,\tsigma_{T,n}^2\ .
\end{equation}

\begin{remark}\rm
Thanks to the simplicity of the Gaussian model, similar estimates may
be obtained for other $\ell^p$-type norms.
\end{remark}
A fundamental aspect of transform coding schemes based on non-linear
approximations such as the one we are describing here is the fact
that the significance maps $\Delta$ have to be encoded together with
the corresponding coefficients. Since the significance map takes
the form of a series of segments of $T$s and segments of $R$s with
various lengths, it is natural to use classical techniques of
run length coding (see for example~\cite{Jayant84digital}, Chapter 10,
for a detailed account) to encode them. The corresponding bit rate
depends crucially on the entropy of the distribution of $T$ and $R$
segments. For the sake of simplicity, let us introduce the entropy
of a binary source with probabilities $(p,1-p)$:
\begin{equation}
h(p) = -p\log_2(p) - (1-p)\log_2(1-p)\ .
\end{equation}
\begin{proposition}
\label{prop:ton.SMrate}
Assume that the initial frequencies of the chains
$X_{\cdot n}$ equal their equilibrium frequencies.
For each frequency bin $n$, the entropy of the distribution of
lengths $L_n$ of $T$ and $R$ segments reads
\begin{equation}
\cH(L_n) = \frac{1-\tpi_n'}{2-\tpi_n-\tpi_n'}h(\tpi_n) +
\frac{1-\tpi_n}{2-\tpi_n-\tpi_n'}h(\tpi_n')\ .
\end{equation}
\end{proposition}
\noindent{\em Proof:} Denote by $L_T$ and $L_R$ the lengths
of $T$ and $R$ segments. From the Markov model $X$ it follows that
$L_T$ and $L_R$ are exponentially distributed:
$$
\Pr{L_T=\ell} = \tpi_n^{\ell-1}(1-\tpi_n)\ ,\quad
\Pr{L_R=\ell} = \tpi_n'^{\ell-1}(1-\tpi_n')\ ,\quad\ell=1,2,\dots
$$
A simple calculation shows that the Shannon entropy of the
random variable $L_T$ is given by
$$
-\sum_{\ell=1}^\infty \Pr{L_T=\ell}\log_2\left(\Pr{L_T=\ell}\right)
= -\tpi_n\log_2(\tpi_n) - (1-\tpi_n)\log_2(1-\tpi_n)=h(\tpi_n)\ ,
$$
and a similar expression for the Shannon entropy of $L_R$.
Now, because of the assumption on the initial frequencies of
the chains $X_{\cdot n}$, and dropping the indices for the
sake of simplicity, we have that
$$
\Pr{X=T}=\frac{1-\tpi'}{2-\tpi-\tpi'}\ ,
$$
and the equality
$$
\cH(L) = \Pr{X=T} \cH(L_T) + \Pr{X=R} \cH(L_R) 
$$
yields the desired result.
\qed

Finally, let us briefly discuss questions regarding the
quantization of coefficients. The simplicity of the model
(Gaussian coefficients, and Markov chain significance map)
makes it possible to obtain elementary rate-distortion estimates.
Indeed, the optimal rate-distortion function for Gaussians
random variables is well known: for a $\cN(0,\sigma^2)$ random
variable,
\begin{equation}
\label{fo:RD.gaussian}
D(R) = \sigma^2\, 2^{-2R}\ .
\end{equation}
Let us assume that the $T$ type coefficients at frequency $n$
are quantized using $R_n$ bits per coefficient.
Using the optimal rate-distortion function~(\ref{fo:RD.gaussian}),
the overall distortion per time frame is given by
$$
D = \sum_{n=0}^{N-1} \frac{\tilde N^{(K)}_n}K\, \tsigma_{T,n}^2\, 2^{-2R_n}\ .
$$
If we are given a global budget of $\overline{R}$ bits
per sample, the optimal bit rate distribution over frequency bins
is obtained by minimizing $\Ex{D}$ with respect to $R_n$,
under the ``global bit budget'' constraint
$$
\Ex{\sum_{n=0}^{N-1} \frac{\tilde N^{(K)}_n}K\, R_n} = N \overline{R} \ ,
$$
the expectation being taken with respect to the significance map $\Delta$.
Assuming for the sake of simplicity that the Markov chain is at equilibrium
(i.e. $\nu_n =\nu_n^{(e)}$ for all $n$),
this yields the following simple expression
\begin{equation}
R_n = \frac{N}{\overline{N}}\overline{R}
+ \frac1{2}\log_2(\tsigma_{T,n}^2) - \frac1{2\overline{N}}
\sum_{m=0}^{N-1} \nu_m^{(e)}\,\log_2(\tsigma_m^2)\ ,
\end{equation}
where we have denoted by
$$
\overline{N} = \sum_{n=0}^{N-1} \nu_n^{(e)}
$$
the average number of $T$ type coefficients per time frame.
As usual in this type of calculation, the so-obtained optimal
value of $R_n$ is generally not an integer number, and an additional
rounding operation is needed in practice. The distortion obtained with
the rounded bit rates is therefore larger than the bound obtained
with the values above. Summarizing this calculation, and plugging
these optimal bit rates into the expression of the distortion, we obtain
\begin{proposition}
\label{th:ton.rate.dist}
With the above notations, the following rate-distortion bound holds:
for a given overall bit budget of $\overline{R}$ bits per $T$ type
coefficient,
\begin{equation}
\Ex{D} \ge \overline{N}\,
\left(\prod_{n=0}^{N-1} \tsigma_n^{2\nu_n^{(e)}}\right)^{1/\overline{N}}\,
2^{-2N\overline{R}/\overline{N}}\ .
\end{equation}
\end{proposition}

\subsection{Parameter and state estimation: algorithmic aspects}
\label{sub:ton.estimation}
Hidden Markov models have been very successful because
there exist naturally associated efficient algorithms for
both parameter estimation and hidden state estimation, respectively
the Expectation-Maximization (EM) and Viterbi algorithms.
However, while these are natural
answers to the estimation problems in general situations, they are not
so natural anymore in a coding setting, as we explain below.

>From a general point of view, an imput signal is first
expanded with respect to an MDCT basis, corresponding to a fixed
time segmentation (segments of approximately 20 msec.) Then, within
larger time frames, the parameters are (re)estimated, as well as the
hidden states. Parameters are refreshed on a regular basis.

\subsubsection{Parameter estimation}
\label{se:ton.param.est}
Given the parameter set $\ttheta$ of the model, the
forward-backward equations allow one to obtain estimates for the
probabilities of hidden states conditional to the observations:
\begin{eqnarray}
\label{fo:pkn.T}
p_{kn}(T)&=&\Pr{X_{k\,n} = T\vert \ttheta,Y_{1:K,n}=y_{1:K,n}}\ ,\\
\label{fo:pkn.R}
p_{kn}(R) &=& \Pr{X_{k\,n}=R\vert \ttheta,Y_{1:K,n}=y_{1:K,n}}
\end{eqnarray}
and the likelihood of the parameters
$$
\cL(\ttheta) = \Pr{Y_{1:K,n}=y_{1:K,n}\vert\ttheta}\ ,
$$
from which new estimates for the parameter set $\ttheta$ may be derived.
\begin{remark}\rm
From a practical point of view, such parameter re-estimation
happens to be quite costly. Therefore, the parameters are generally
re-estimated on a larger time scale, taking several consecutive
windows into account.
\end{remark}
\begin{remark}\rm
\label{rem:freq.decay}
For practical purpose, it is generally more suitable to
restrict the parameter set $\ttheta$ to a smaller subset.
The following two assumptions proved to be quite adapted
to the case of audio signals:
\begin{itemize}
\item[{\em i.}] The variances may be assumed to be multiple of
a single reference value, implementing some ``natural'' decay of
MDCT coefficients with respect to frequency. For example, we
generally used expressions of the form
$$
\tsigma_{s,n} = \frac{\tsigma_s}{1 + (\frac{n}{n_0})^\alpha} ,\quad s=T,R
$$
$n_0\in\RR^+$ being some reference frequency bin, $\sigma_s$ a reference
standard deviation for state $s$ and $\alpha \in\RR^+$ some constant
controlling the decay (typical value being $\alpha=1$).
Without such an assumption,
frequency bins are completely independent of each other, and the estimation
algorithm generally yields $T$ type coefficients in all bins,
which is not realistic,
\item[{\em ii.}] For each frequency bin, the initial frequencies
$\nu_n$ of the considered Markov chain are generally assumed to
equal the equilibrium frequencies $\nu_n^{(e)}$.
\end{itemize}
\end{remark}

\subsubsection{State estimation}
Viterbi's algorithm is generally considered the natural answer to
the state estimation problem. It is a dynamic programing
algorithm, which yields Maximum a posteriori (MAP) estimates
$$
\hat x_{1:K\, n} = \arg\max \Pr{X_{1:K\,n}
= x_{1:K\,n}\vert y_{1:K\,n},\ttheta}\ ,
$$
for each frequency bin $n$. However, the number of so-obtained
coefficients in a given state ($T$ or $R$) cannot be controlled
a priori when such an algorithm is used, which turns out to be
a severe limitation in a signal coding perspective. In addition,
Viterbi's algorithm requires that accurate estimates of the model's
parameters are available, which will not necessarily be the case
if the parameter estimates are refreshed on a coarse time scale
(see above.)

Therefore, we also consider, as an alternative to Viterbi's algorithm,
an {\em a posteriori probabilities thresholding} method, which is
computationally far simpler, and allows a fine rate control.
More precisely, given a prescribed rate $N_{ton}$,
%the selected $T$-type coefficients are those whose a posteriori probability
%$p_{kn}(T)$ in~(\ref{fo:pkn.T}) are the $N_{ton}$ largest ones:
\begin{itemize}
\item[{\em i.}] Sort the MDCT coefficients
$y_{kn}=\langle x,w_{kn}\rangle$ in order of decreasing a posteriori
probability $p_{kn}(T)$ in~(\ref{fo:pkn.T}),
\item[{\em ii.}] Keep the $N_{ton}$ first sorted coefficients.
\end{itemize}
In this way, for an average bit rate $\overline{R}$ and a
prescribed ``tonal'' bit budget, a number $N_{ton}$
of MDCT coefficients to be retained may be estimated, and
the $N_{ton}$ coefficients with largest a posteriori probability
are selected.

\subsection{Numerical simulations}

\begin{figure}
\centerline{
\includegraphics[height=11cm,width=2.5cm,angle=-90]{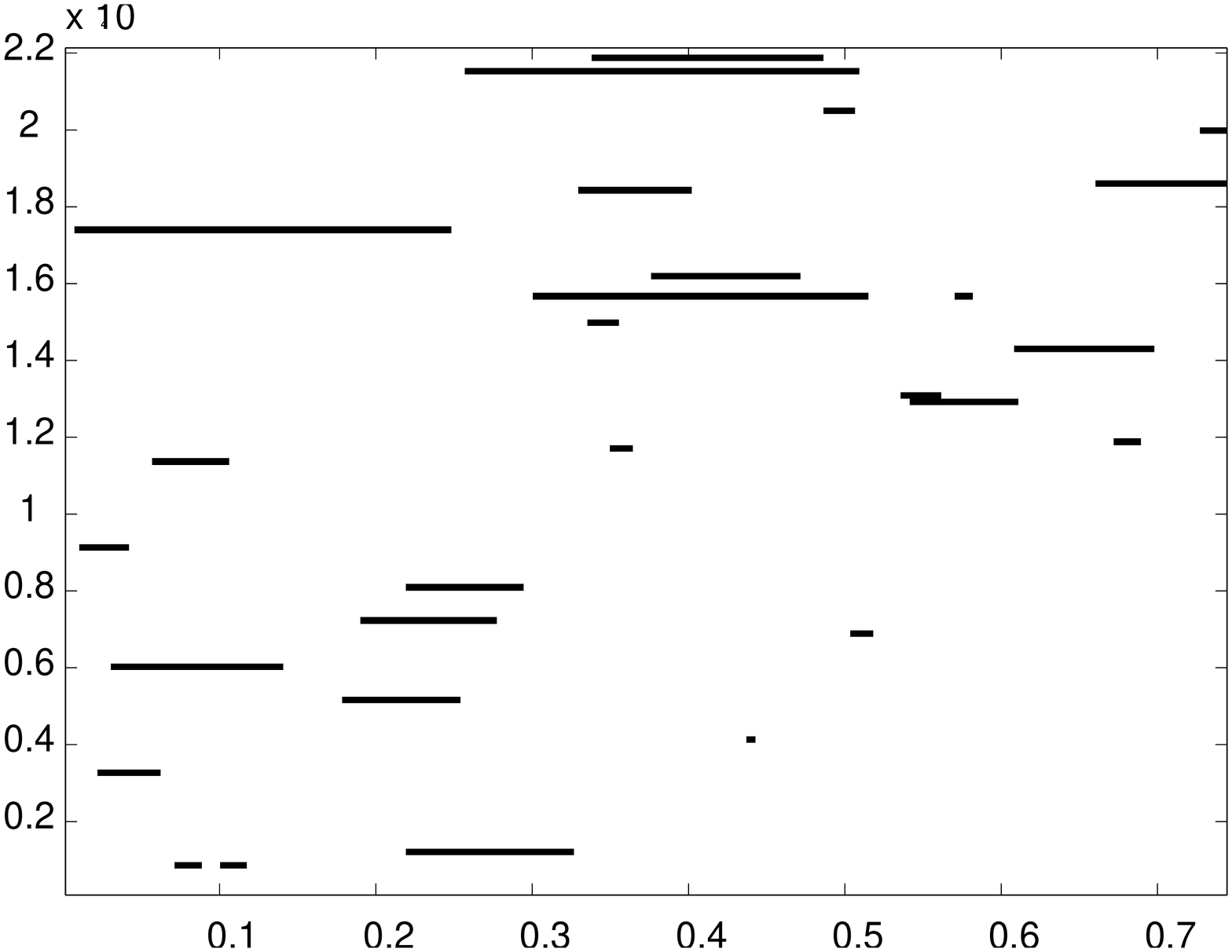}
}
\centerline{
\includegraphics[height=11cm,width=2.5cm,angle=-90]{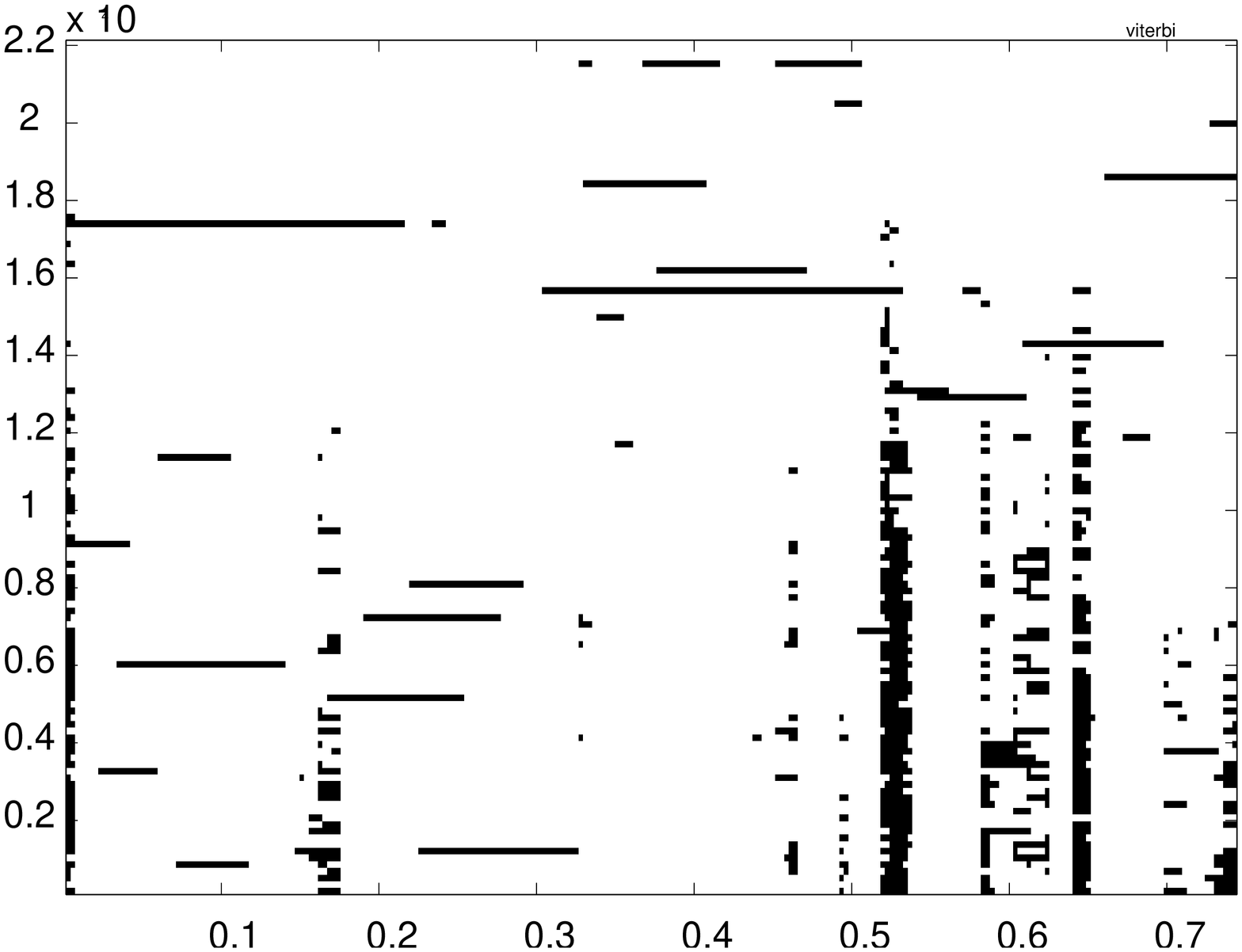}
}
\smallskip
\centerline{
\hskip2mm
\includegraphics[height=10.8cm,width=6cm,angle=-90]{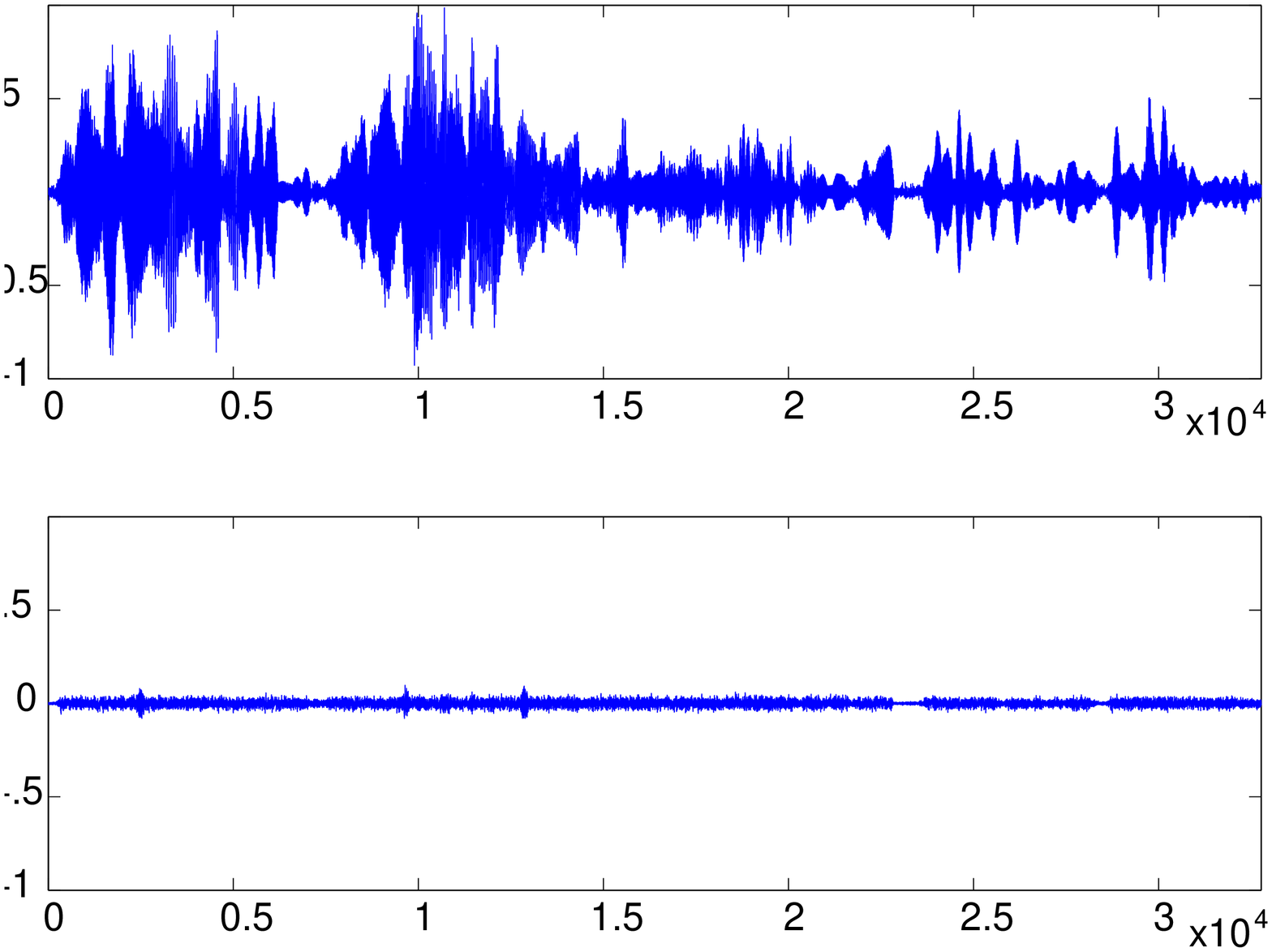}
}
\caption{Estimating a tonal layer from simulated signal;
from top to bottom:
simulated significance map, estimated significance map
(estimation via the Viterbi algorithm),
estimated tonal signal, estimated residual signal.}
\label{fi:sim.ton.est}
\end{figure}

As a first test of the model and the estimation algorithms,
we generated realizations of the structured harmonic mixture
of Gaussians model described above, and used the corresponding
estimation algorithms. We simulated a signal according to the
``tonal + residual'' Markov model as above, with about 3.1\%
$T$-type coefficients. We show in {\sc Figure}~\ref{fi:sim.ton.est}
the result of the estimation of the tonal layer using EM parameter
estimation, and state estimation via the Viterbi algorithm.
As may be seen, the significance map
is fairly well estimated, except in regions where the signal has
little energy, which was to be expected. In these regions, the
algorithm detects spurious (vertical) tonal structures,
which results in an increase of the percentage of
$T$ type coefficients (about 4.1\% instead of 3.2\% for that example.)
However, since this effect appears only in regions where the
signal has small energy, this does not affect tremendously
the estimated signal, which is very close to the simulated one
(not shown here.)

For the sake of comparison, we display in {\sc Figures}~\ref{fi:sim.ton.est2}
and~\ref{fi:sim.ton.est3} some examples of tonal layer estimation
using the thresholding algorithm instead of the Viterbi algorithm, for
various values of the threshold. The simulation presented in
{\sc Figure}~\ref{fi:sim.ton.est2} corresponds to 1\% retained coefficients, while
the simulation presented in {\sc Figure}~\ref{fi:sim.ton.est3} corresponds
to 3\% retained coefficients. As expected, the significance map in
{\sc Figure}~\ref{fi:sim.ton.est2} appears much terser than the ``true'' one,
while the one in {\sc Figure}~\ref{fi:sim.ton.est3} is much closer (percentage
of retained coefficients significantly larger than the true one
yield spurious tonal structures.) This results in tonal components
which were not correctly captured, and appear in the residual signal
of {\sc Figure}~\ref{fi:sim.ton.est2}. This is not the case any more
when the threshold is set to a more ``realistic'' value, as may be seen
in the tonal and residual layers of {\sc Figure}~\ref{fi:sim.ton.est3}.
In that case, the residual only features a small spurious component.
Notice that even though significantly less coefficients
are retained, the overall shape of the estimated signal is quite
good.

\begin{figure}
\centerline{
\includegraphics[height=11cm,width=2.5cm,angle=-90]{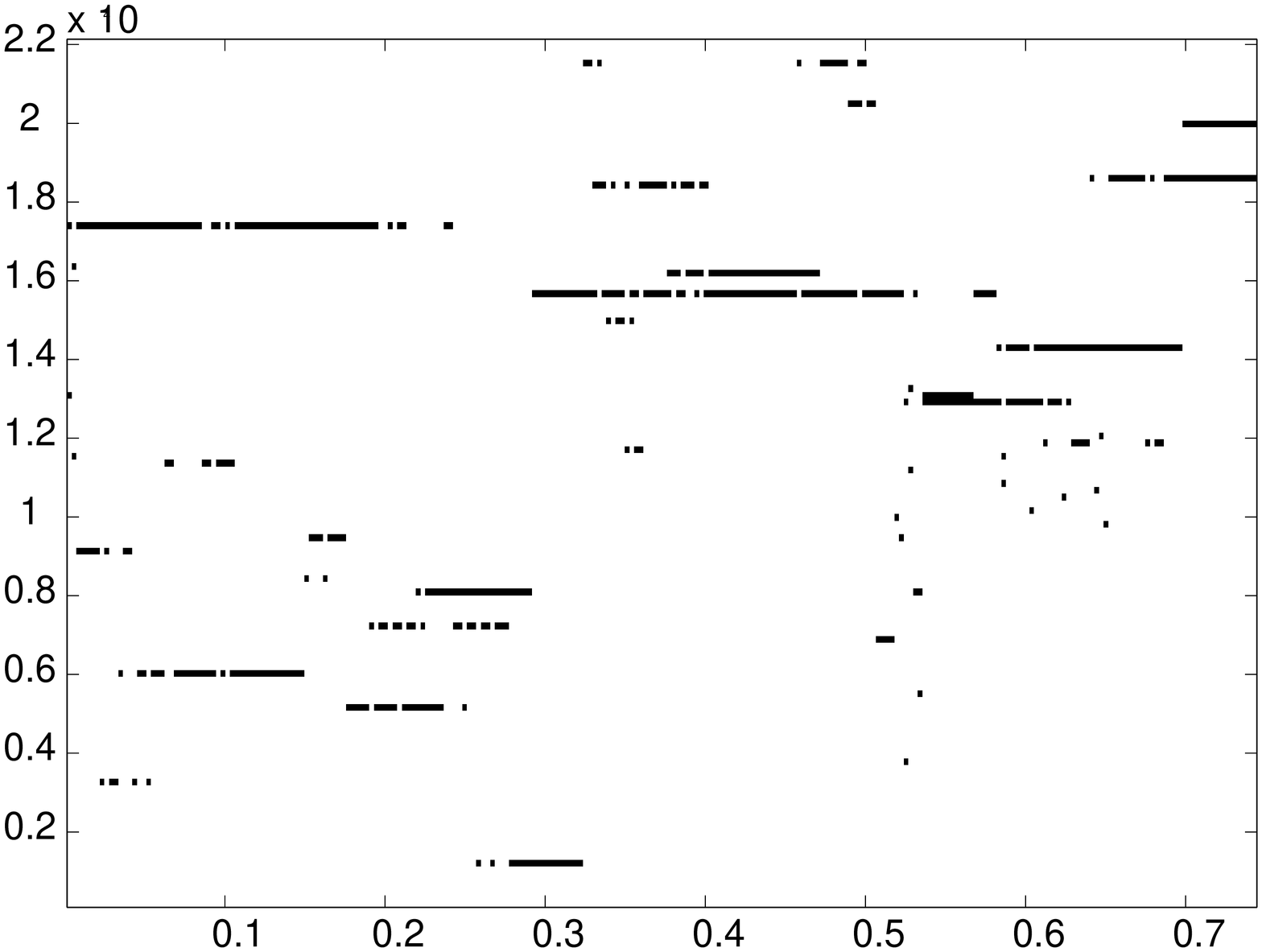}
}
\centerline{
\hskip2mm
\includegraphics[height=10.8cm,width=6cm,angle=-90]{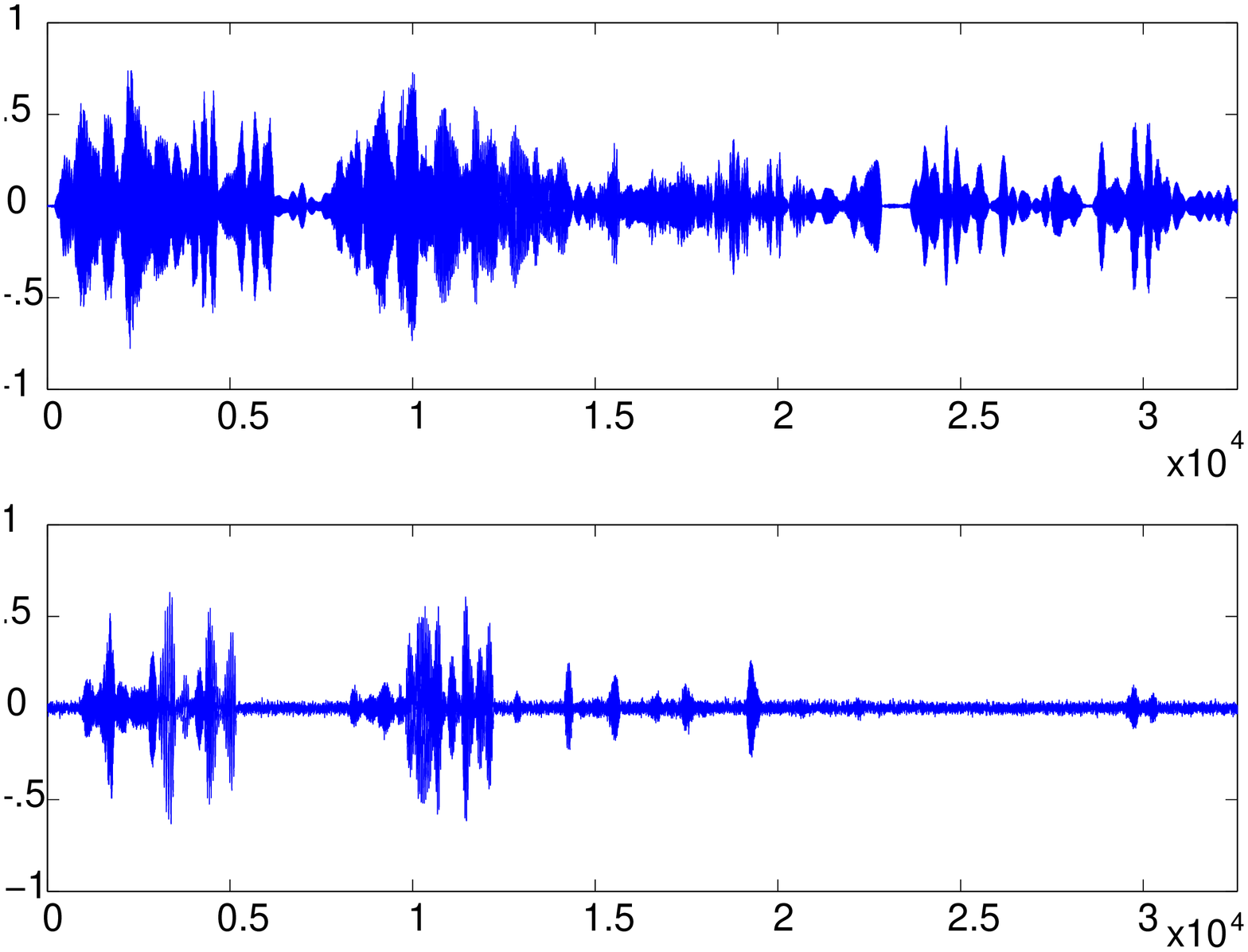}
}
\caption{Estimating a tonal layer from simulated signal (same as in
{\sc Figure}~\ref{fi:sim.ton.est});
from top to bottom: estimated significance map
(estimated via the posterior probability thresholding algorithm,
using 1\% coefficients);
estimated tonal signal, estimated residual signal.}
\label{fi:sim.ton.est2}
\end{figure}
\begin{figure}
\centerline{
\includegraphics[height=11cm,width=2.5cm,angle=-90]{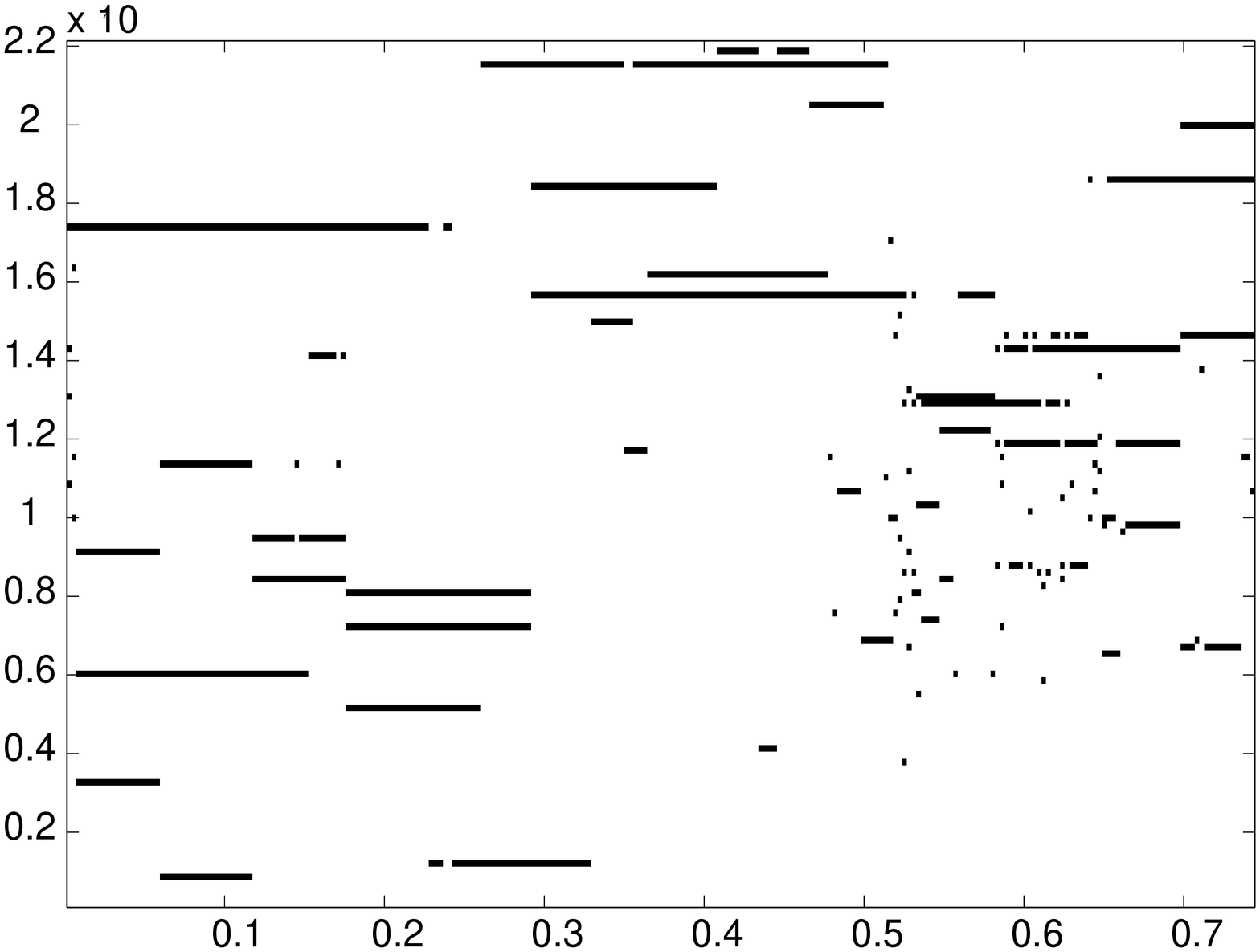}
}

\medskip
\centerline{
\hskip2mm
\includegraphics[height=10.8cm,width=6cm,angle=-90]{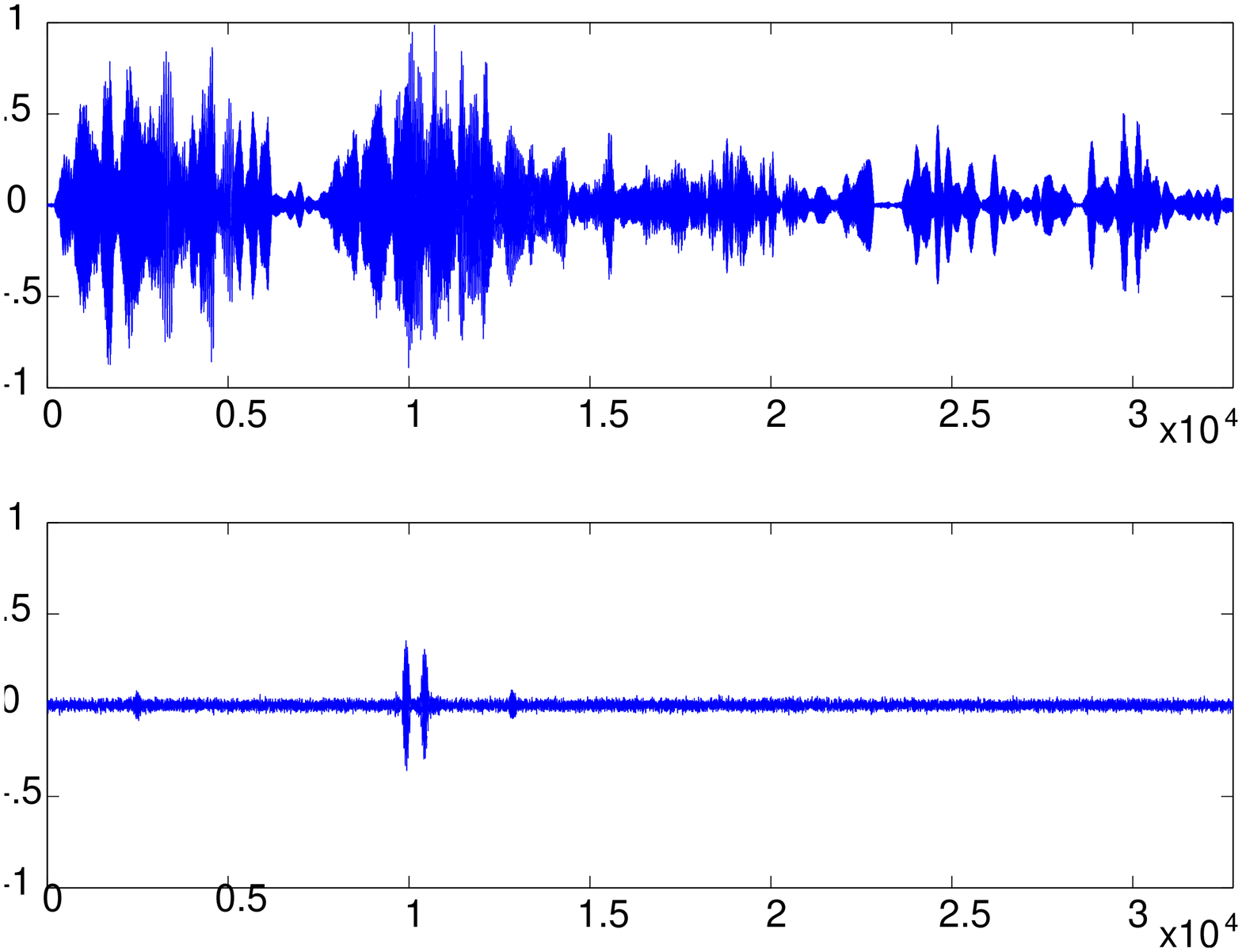}
}
\caption{Estimating a tonal layer from simulated signal
(same as in {\sc Figure}~\ref{fi:sim.ton.est});
from top to bottom:
estimated significance map
(estimated via the posterior probability thresholding algorithm,
using 3\% coefficients);
estimated tonal signal, estimated residual signal.}
\label{fi:sim.ton.est3}
\end{figure}

\begin{remark}
\label{rem:viterbvspost}\rm
Clearly, the posterior probability thresholding method only provides
an approximation of the ``true'' tonal layer (which is provided by
the Viterbi algorithm), whose precision depends on the choice of the
threshold, i.e. the bit rate allocated to the tonal layer.
Controlling the relation between the bit rate and the precision of
the approximation would lead to a rate distortion theory for
the ``functional'' part of the tonal coder. Such a theory seems
extremely difficult to develop, and so far we could only study it
by numerical simulations (not shown here.) 
\end{remark}
%
%%%%
%
\section{Structured Markov model for transient}
\label{se:transient}
\subsection{Hidden wavelet Markov tree model}
We now turn to the description of the transient model, which was
partly presented in~\cite{Molla03hidden}. The latter
exploits the fact that wavelet bases are ``well adapted''
for describing transients, in the sense that these
generally yield scale-persistent chains of significant
wavelet coefficients. We start from a
multiresolution analysis (see for
example~\cite{Mallat98wavelet,Vetterli95wavelets})
and the corresponding wavelet $\psi\in\Ltwo$, scaling
function $\phi\in\Ltwo$ and wavelet basis, defined by
$$
\psi_{jk}(t) = 2^{-j/2}\,\psi\!\left(2^{-j}t-k\right)\ ,\quad
j,k\in\ZZ\ .
$$
Given $x\in\Ltwo$, its wavelet coefficients
$d_{jk}=\langle x,\psi_{jk}\rangle$ are naturally labelled
by a dyadic tree, as in {\sc Figure}~\ref{fi:dyadtree}, in
which it clearly appears that a given wavelet coefficient
$d_{jk}$ may be given a pair of children
$d_{j+1\ 2k}$ and $d_{j+1\ 2k+1}$.
For the sake of simplicity, we shall sometimes
collect the two indices $j,k$
into the scale-time index $\lambda=(j,k)$.

For the sake of simplicity, we consider a fixed time interval,
and a signal model involving finitely many scales, of the form
\begin{equation}
\label{fo:nton.model}
x = S_{J0}\phi_{J0} + \sum_{j=1}^J \sum_{k=0}^{2^{J-j}-1} D_{jk}\psi_{jk}\ ,
\end{equation}
involving
$$
N(J) = 2^J -1
$$
random wavelet coefficients\footnote{The scaling function
coefficients $S_{J0}$ are generally irrelevant for audio signals, and do
not deserve much modelling effort.}, whose distribution is
a gaussian mixture governed by a hidden random variable.
%the latter has two possible states: $T$ for ``transient'' and
%$R$ for ``residual''. Transitions from a state to another are
%governed by a (hidden) Markov chain.

\begin{figure}
\centerline{
\psfig{file=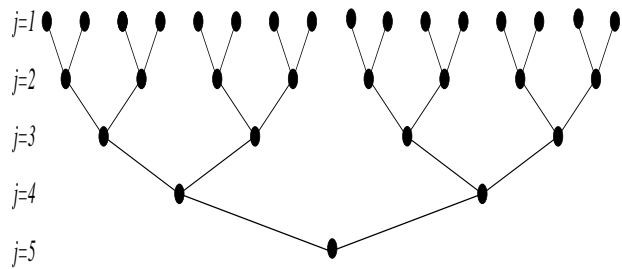,height=35mm,width=8cm}
}
\caption{Wavelet coefficients tree.}
\label{fi:dyadtree}
\end{figure}

More precisely, distribution of the wavelet coefficients
$D_{jk}$ depends on a hidden state $X_{jk}\in\{T,R\}$
($T$ stands for ``transient'', and $R$ for ``residual''.)
At each scale $j$, the $T$-type
coefficients are modelled by a centered normal distribution
with (large) variance $\sigma_{T,j}^2$.
The $R$-type coefficients are modelled by a centered normal
distribution with (small) variance $\sigma_{R,j}^2$.

The distribution of hidden states is given by a ``coarse to fine''
Markov chain, characterized by a $2\times 2$ transition matrix,
and the distribution of the coarsest scale state.
In order to retain only connected trees, we
impose a {\em taboo transition}: the transition $R\to T$ is forbidden.
Therefore, the transition matrix assumes the form
$$
P_j = \left(
\begin{array}{cc}
\pi_j&1-\pi_j\\
0&1
\end{array}
\right)
$$
where $\pi_j$ denotes the {\em scale persistence} probability,
namely the probability of transition $T\to T$ at scale $j$:
$$
\pi_j =\Pr{X_{j-1,\ell}=T\vert X_{j,k}=T}\ ,\ell=2k,2k+1\ .
$$
The hidden Markov process is completely determined by the
set of matrices $P_j$ and the ``initial'' probability distribution,
namely the probabilities $\nu=(\nu_T,\nu_R)$ of states
at the maximum considered scale $j=J$.
The complete model is therefore characterized by the numbers
$\pi_j$, $\nu$, and the emission probability densities:
$$
\rho_S(d) = \rho(d\vert X=S)\ ,\quad S=T,R\ .
$$
In the sequel, we shall always assume that the persistence
probabilities are scale independent\footnote{This is actually
quite a strong assumption, which has the advantage of reducing the
number of parameters to estimate. Alternative choices can also
be considered, for example controlling the growth of the number of
significant coefficients across scales by setting $\pi_i = \pi_0^i$ for
some constant $\pi_0$.}:
$$
\pi_i =\pi\ ,\quad\forall i\ .
$$
According to our choice (centered Gaussian distributions),
the latter are completely characterized by their variances
$\sigma_{T,j}^2$ and $\sigma_{R,j}^2$.
All together, the model is completely specified by
the parameter set
\begin{equation}
\label{fo:tr.param}
\theta = \{\nu ,\, \pi ,\, \sigma_{T,j} ,\, \sigma_{R,j} ,\; j~=~1\dots J \}\ ,
\end{equation}
which leads to the definition of {\em transient significance map}
(termed transient feature in~\cite{Molla03hidden})
\begin{definition}
\label{def:transient}
Let the parameter set in~(\ref{fo:tr.param}) be fixed, and let $x$
denote a signal given by a hidden Markov tree model as
in~(\ref{fo:nton.model}) above. Consider the random set
\begin{equation}
\Lambda = \left\{(j,k),j=1,\dots j,\ k=0,\dots 2^j-1\vert X_{jk}=T\right\}
\ .
\end{equation}
$\Lambda$ is called the {\em transient significance map} of $x$.
The corresponding transient layer of $x$ is defined as
\begin{equation}
x_{tr} = \sum_{(j,k)\in\Lambda} D_{jk} \psi_{jk}\ .
\end{equation}
\end{definition}
From this definition, one may easily derive estimates on various coding
rates. The key point is the following immediate remark. Let $N_j$ denote
the number of $T$-type coefficients at scale $j$, and let
$$
N = \sum_{j=1}^J N_j
$$
the total number of $T$-type coefficients at scale $j$.
The following result is fairly classical in branching processes theory
(see for example~\cite{Doob53stochastic,Karlin97first}.)
\begin{proposition}
Let $x$ denote a signal given by a hidden Markov tree model as
in~(\ref{fo:nton.model}) above. Then the number $N$ of
$T$ type coefficients is given by a Galton-Watson process. In
particular, one has
\begin{equation}
\Ex{N_j} = \nu (2\pi)^{J-j}\ ,\quad
\overline{N} := \Ex{N} = \nu\,\frac{(2\pi)^J -1}{2\pi -1}
\end{equation}
(with the obvious modification for the case $\pi=1/2$.)
\end{proposition}
Therefore, it is obvious to obtain estimates for the
energy of a transient layer:
\begin{corollary}
\label{prop:tr.SMrate}
The average energy of the transient layer of a signal $x$ reads
\begin{equation}
\Ex{\sum_{j,k;X_{jk}=T} |D_{jk}|^2} = 
\nu \sum_{j=1}^{j=J} \sigma_j^2 (2 \pi)^{J-j}\ .
\end{equation}
\end{corollary}
Another simple consequence is the following a priori estimate for the
cost of significance map encoding. It is known that it is possible to encode
a binary tree at a cost which is linear in the number of nodes.
We use the following strategy for encoding the tree $\Lambda$ (even though
it is not optimal, it has the advantage of being simple. Improvements may
be obtained by using entropy coding techniques, taking advantage of the
probability distribution of trees, which is known as soon as the persistence
probability $\pi$ is known.)
We associate with each node of $\Lambda$ a pair of bits,
set to 0 or 1 depending on whether the left and right children of
the node belong to $\Lambda$ or not. Therefore, $R_{SM}$
is not larger than twice the number of nodes of $\Lambda$,
i.e. the number of $T$-type coefficients.
Therefore, we immediately deduce
\begin{corollary}
\label{th:tree_encoding}
Given the set of parameters $\theta$, and the corresponding
Hidden Markov wavelet tree model,
let $R_{SM}$ denote the number of bits necessary to encode the
significance map of a transient wavelet coefficients tree,
as above. Then we have
$$
\Ex{R_{SM}} \le \left\{ \begin{array}{ll}\dd
         2\nu \times \frac{1-(2 \pi)^J}{1-2 \pi} &\hbox{ if }\pi \neq 0.5,\\
         2\nu J& \hbox{ if }\pi = 0.5 .
         \end{array}\right.
$$
\end{corollary}

\medskip

The simplicity of the transient model (i.e. Galton-Watson
significance map, and Gaussian $T$ coefficients) makes it possible to
derive simple rate-distortion estimates, along lines similar to the
ones we followed for the tonal layer.
Assume that the $T$ type coefficients at scale $j$
are quantized using $R_j$ bits. Assuming~(\ref{fo:RD.gaussian}),
the overall distortion is given by
$$
D = \sum_{j=1}^J N_j \sigma_j^2 2^{-2R_j}\ .
$$
Suppose we are given a global budget of $\overline{R}$ bits
per sample.
Minimizing $\Ex{D}$ with respect to $R_j$, under the ``global bit budget''
constraint
$$
\Ex{\sum_{j=1}^J N_j R_j} = N(J) \overline{R} 
$$
yields the following simple expression
\begin{equation}
R_j = \frac{N(J)}{\overline{N}}\overline{R}
+ \frac1{2}\log_2(\sigma_j^2) - \frac1{2}
\frac{2\pi -1}{(2\pi)^J-1} \sum_{j=1}^J (2\pi)^{J-j}\,\log_2(\sigma_j^2)\ .
\end{equation}
Therefore, plugging this expression into the optimal
rate-distortion function~(\ref{fo:RD.gaussian}),
we obtain the following rate-distortion estimate
\begin{proposition}
\label{th:rate.dist}
With the same notations as before, we have the following estimate:
for a given overall bit budget of $\overline{R}$ bits per $T$ type
coefficient, the distortion is such that
\begin{equation}
\Ex{D} \ge \overline{N}\,
\left(\prod_{j=1}^J \sigma_j^{2\overline{N_j}}\right)^{1/\overline{N}}\,
2^{-2N(J)\overline{R}/\overline{N}}\ ,
\end{equation}
where we have set
$$
\overline{N_j} = \nu\,(2\pi)^{J-j}\ ,\quad
\overline{N} =\nu\,\frac{(2\pi)^J-1}{2\pi-1}\ .
$$
\end{proposition}

\subsection{Parameters and state estimation}
As in the case of the tonal layer, the parameter estimation and
the hidden state estimation may be realized through standard EM
and Viterbi type algorithms. These algorithms are mainly
based upon adapted versions of the above mentioned forward-backward
algorithm: the so-called ``upward-downward'' algorithm, proposed by
Crouse and collaborators in~\cite{Crouse98wavelet}. Actually,
we rather used a variant, the downward-upward algorithm, due
to Durand and Gon\c calves~\cite{Durand01statistical}, which
provides a better control of numerical accuracy of the computations.
As a result, the algorithm provides estimates for quantities such
as the hidden states probabilities
$$
\Pr{X_{jk}=s\big| D_{1:2^J-1}=d_{1:2^J-1}1,\theta}
$$
and the likelihood
$$
\cL = \rho_{D_{1:2^J-1}}\left(d_{1:2^J-1}\big|
X_{{1:2^J-1}},\theta\right)\ .
$$
\subsubsection{Parameters estimation}
\label{se:tr.param.est}
The parameter estimation goes along lines similar to the
ones outlined in Section~\ref{se:ton.param.est}
(see also~\cite{Molla03hidden} for additional details.) Again, since
the parameter estimation procedure, involving upward-downward
algorithm, is quite costly, it is done simultaneously on several
consecutive time windows (i.e. several consecutive trees), and
parameters are ``refreshed'' on larger time scales.

\subsubsection{Hidden states estimation}
Again, the situation is very similar to the situation encountered
when dealing with the tonal layer. The ``Viterbi-type'' algorithm described
in~\cite{Durand01statistical} theoretically provides an estimate for
``the'' transient significance map, and therefore the transient layer.
However, it does not allow one to control the number of selected
coefficients (the rate), and is therefore not appropriate in a context
of variable bit rate coder. Hence, we rather turn to the (also
computationally simpler) alternative, using thresholding of a posteriori
probabilities.

The upward-downward algorithm provides estimates for the probabilities
$$
p_{jk}(T) = \Pr{X_{jk}=T\big|D_{1:2^J-1}=d_{1:2^J-1},\theta}\ .
$$
Therefore, the corresponding tree nodes may be sorted according
to the latter (in decreasing order.) For a given transient bit budget,
a maximal number of nodes to be retained $N_{tr}$ may be estimated, and
the nodes with largest ``transientness'' probability $p_{jk}(T)$
are selected, and the corresponding transient layer is reconstructed.

\subsection{Numerical simulations}
As for the case of the tonal layer, it is easy to perform
numerical simulations of the model to evaluate the performances of the
estimation algorithms. We display in {\sc Figure}~\ref{fi:sim.trans.est}
the results of such simulations, using EM algorithm for parameter
estimation, and the Viterbi algorithm for hidden states estimation.
As may be seen from the plots, the significance tree
and the transient layer are quite well estimated.

\begin{figure}
\centerline{
\includegraphics[height=11cm,width=11cm,angle=-90]{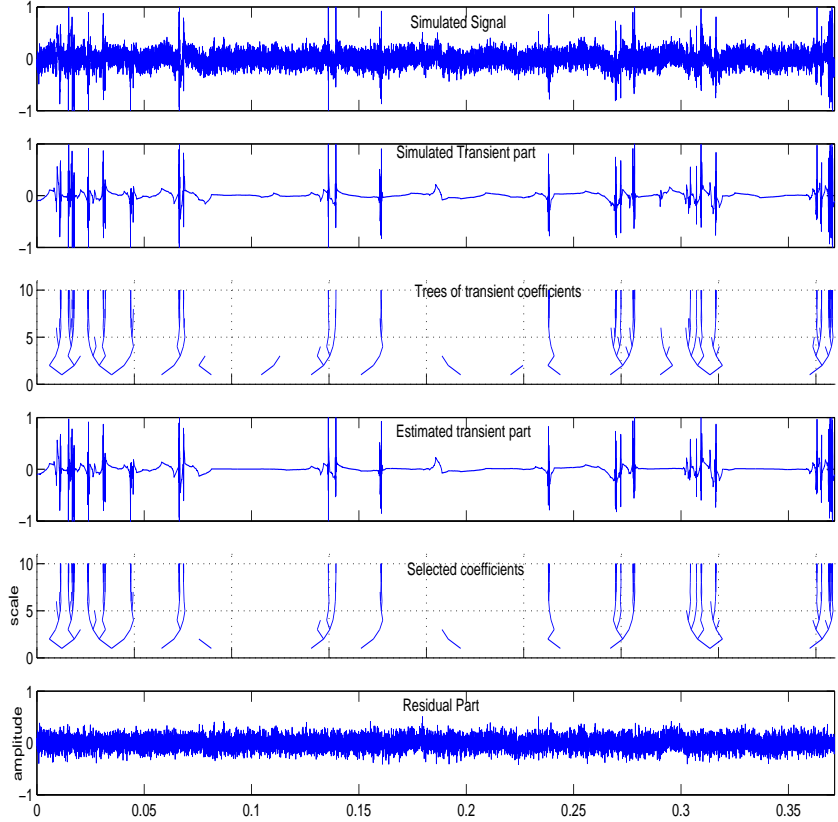}
}
\caption{Estimating a transient layer from simulated signal;
from top to bottom: simulated signal, simulated transient layer,
simulated significance tree, estimated transient layer
(estimation via the Viterbi algorithm),
estimated significance tree,
estimated residual signal.}
\label{fi:sim.trans.est}
\end{figure}

Again, using the posterior probability thresholding method instead
of the Viterbi method yields approximate transient layer, and the
discussion of Remark~\ref{rem:viterbvspost} still hold true.

%\subsubsection{Miscellaneous comments}

\section{The ``tonal \lowercase{vs} transient'' balance}
\label{se:balance}
We have described in Sections~\ref{se:tonal} and~\ref{se:transient}
two models for tonal and transient layers in audio signals, and
corresponding estimation algorithms. One of the main aspects of the
latter is that the hidden states estimation is based on thresholding
of a posteriori probabilities rather than on a global Viterbi-type
estimation, which allows to accomodate any bit rate prescribed
in advance.

However, as stressed in the
Introduction, and described in more detail in the subsequent section,
we develop a coding approach based upon recursive estimations of
tonal and transient layers. We describe below an approach for
pre-estimating the relative sizes of the tonal and transient
layers, in order to balance the bit budget between the two layers
prior to estimation. The reader interested in more details is invited to
refer to~\cite{Molla03determining}.

\subsection{Pre-estimating the ``sizes'' of the tonal and transient layers}
Consider a signal assumed for simplicity to be of the
form~(\ref{fo:sig.model}), with unknown values of $|\Delta|$
and $|\Lambda|$, we seek estimates for the ``transientness''
and ``tonality'' indices
\begin{equation}
\label{fo:indices}
I_{ton} = \frac{|\Delta|}{|\Delta|+|\Lambda|}\ ;\qquad
I_{tr} = \frac{|\Lambda|}{|\Delta|+|\Lambda|}\ ,
\end{equation}
or alternatively, the proportion of the signal's energy contained
in the tonal and transient layers.
For simplicity, we limit ourselves to the finite dimensional situation,
and propose a procedure very much in the spirit of the information
theoretic approaches advocated by M.V. Wickerhauser and
collaborators~\cite{Trgo95relation,Wickerhauser94adapted}.
\begin{definition}
Let $\cB = \{e_n,n\in S\}$ be an orthonormal basis
of a given $N$-dimensional signal space $\cE$.
The logarithmic dimension of $x\in\cE$ in the basis $\cB$
is defined by
\begin{equation}
\cD_\cB(x) = \frac1{N}\,\sum_{n\in S}
\log_2\left(|\langle x,e_n\rangle|^2\ \right)
\end{equation}
\end{definition}
We aim to show that such quantity may provide the desired
estimates, under suitable assumptions on the signal (sparsity) and
the considered bases (incoherence.)
Elementary calculations show that in the framework of
the signal models~(\ref{fo:sig.model}), one has the following
\begin{lemma}
Given an orthonormal basis $\cB = \{e_n,n\in S\}$,
assuming that the coefficients $\langle x,e_n\rangle$
of $x\in\cE$ are $\cN(0,\sigma_n^2)$ random variables, one has
\begin{equation}
\Ex{\cD_\cB(x)} = %\left(1+\frac\gamma{\log(2)}\right) 
C+ \frac1{N}\sum_{n\in S}
\log_2(\sigma_n^2)
\end{equation}
where $C= 1+\gamma/\ln(2)$ ($\gamma\approx .5772156649$ being
Euler's constant.)
\end{lemma}

Consider now the model~(\ref{fo:sig.model}), and assume that
the coefficients $\alpha_\lambda,\lambda\in\Lambda$
and $\beta_\delta,\delta\in\Delta$ are respectively $\cN(0,\sigma_\lambda^2)$
and $\cN(0,\tsigma_\delta^2)$ independent random variables.
Then the coefficients
$$
a_\lambda = \langle x,\psi_\lambda\rangle\ ;\quad
b_\delta = \langle x,w_\delta\rangle\ ,
$$
are centered normal random variables, whose variances depends
on whether $\lambda\in\Lambda$ (or $\delta\in\Delta$) or not.
For example, in the case of the $a_\lambda$ coefficients,
\begin{equation}
\hbox{var}\{a_\lambda\} = \left\{
\begin{array}{ll}
\sigma_\lambda^2 + \sum_{\delta \in\Delta}\tsigma_\delta^2
|\langle x,w_\delta\rangle|^2 &\hbox{ if }\lambda\in\Lambda\\
\sum_{\delta \in\Delta}\tsigma_\delta^2
|\langle x,w_\delta\rangle|^2 &\hbox{ if }\lambda\not\in\Lambda\ ,\\
\end{array}
\right.
\end{equation}
which yields
\begin{equation}
\label{fo:log.dim.1}
\Ex{\cD_\Psi(x)} = C + %\left(1+\frac\gamma{\log(2)}\right) +
\frac1{N}\log_2\!\!
\left(\prod_{\lambda\in\Lambda}\!\!\left(\sigma_\lambda^2 +\!\!
\sum_{\delta\in\Delta}
\tsigma_\delta^2 |\langle \psi_\lambda,w_\delta\rangle|^2\right)
\prod_{\lambda'\not\in\Lambda}\!\!\left(\sum_{\delta\in\Delta} \tsigma_\delta^2
|\langle \psi_{\lambda'},w_\delta\rangle|^2\right)\right)\ ,
\end{equation}
and a similar expression for the logarithmic dimension
$\cD_W(x)$ with respect to the $W=\{w_\delta\}$ basis.

For the sake of simplicity, we now assume that $\sigma_\lambda =\sigma$,
$\forall\lambda\in\Lambda$ and $\tsigma_\delta =\tsigma$,
$\forall\delta\in\Delta$. Introduce the Parseval weights
\begin{equation}
p_\lambda(\Delta) = \sum_{\delta\in\Delta}
|\langle w_\delta,\psi_{\lambda}\rangle|^2\ ,
\quad
\tp_\delta(\Lambda) = \sum_{\lambda\in\Lambda}
|\langle w_\delta,\psi_{\lambda}\rangle|^2\ .
\end{equation}
The Parseval weights provide information regarding the
``dissimilarity'' of the two considered bases.
The following property is a direct consequence of Parseval's formula:
\begin{lemma}
\label{le:parseval}
With the above notations, the Parseval weights satisfy
$$
0\le p_\lambda(\Delta)\le 1\ ,\quad
0\le \tp_\delta(\Lambda)\le 1\ .
$$
\end{lemma}
Introduce the {\em relative redundancies} of the bases
$\Psi$ and $W$ with respect to the significance maps
\begin{equation}
\epsilon(\Delta) = \max_{\lambda\in\Lambda} p_\lambda(\Delta)\ ,
\quad
\tepsilon(\Lambda) = \max_{\delta\in\Delta} p_\delta(\Lambda)\ .
\end{equation}
These quantities carry information similar to the one carried
by the Babel function used in~\cite{Tropp02greed} for example.
One then obtains simple estimates for the logarithmic
dimension~\cite{Molla03determining}.
\begin{proposition}
With the above notations, assuming that the significant coefficients
$\alpha_\lambda,\lambda\in\Lambda$ and
$\beta_\delta,\delta\in\Delta$ are {\em i.i.d.}
$\cN(0,\sigma^2)$ and $\cN(0,\tsigma^2)$ normal variables respectively,
one has the following bound
\begin{eqnarray}
\label{fo:sup.bound}
\quad\Ex{\cD_\Psi(x)} &\ge& C + \frac{|\Lambda|}N\log_2(\sigma^2) +
\log_2\left(\prod_{\lambda'\not\in\Lambda}\! \left(
\tsigma^2 p_{\lambda'}(\Delta)\right)^{1/N}\right)\\
\label{fo:inf.bound}
\quad\Ex{\cD_\Psi(x)} &\le&C + \frac{|\Lambda|}N
\log_2(\sigma^2 +\epsilon(\Delta)\tsigma^2) +
\log_2\left(\prod_{\lambda'\not\in\Lambda}\! \left(
\tsigma^2 p_{\lambda'}(\Delta)\right)^{1/N}\right)\ .
\end{eqnarray}
Exchanging the roles of $\Delta$ and $\Lambda$, 
a similar bound is obtained for $\cD_W(x)$.
\end{proposition}
%\n{\em Proof: } the proposition follows directly from the fact that
%in such a situation, equation~(\ref{fo:log.dim.1}) reduces to
%\begin{equation}
%\label{fo:log.dim.2}
%\Ex{\cD_\Psi(x)} = C + 
%\log_2\left(\prod_{\lambda\in\Lambda}\! \left(\sigma^2 +
%\tsigma^2 p_\lambda(\Delta)\right)^{1/N}
%\prod_{\lambda'\not\in\Lambda}\! \left(
%\tsigma^2 p_{\lambda'}(\Delta)\right)^{1/N}\right)\ ,
%\end{equation}
%from Lemma~\ref{le:parseval} and the definition of
%$\epsilon(\Delta)$.\foorp

\medskip
\n At this point, several comments have to be made.
\begin{itemize}
\item[{\em a.}] The bounds in Equations~(\ref{fo:sup.bound})
and~(\ref{fo:inf.bound}) differ by
$|\Lambda|\log_2(1+\epsilon(\Delta)\tsigma^2/\sigma^2)/N$.
Let us temporarily assume that this term may be neglected
(see comment {\em b.} below for more details.)
The behavior of $\Ex{\cD_\Psi(x)}$
is therefore essentially controlled by
$$
\log_2\left(\prod_{\lambda'\not\in\Lambda}\! \left(
\tsigma^2 p_{\lambda'}(\Delta)\right)^{1/N}\right)
$$
Such an expression is not easily understood, but a first
idea may be obtained by replacing $p_{\lambda'}(\Delta)$ by its
``ensemble average''
$$
\frac1{N}\sum_{\lambda=1}^N p_{\lambda}(\Delta) = 
\frac1{N}\sum_{\lambda=1}^N\sum_{\delta\in\Delta}
|\langle w_\delta,\psi_\lambda\rangle|^2 =
\frac1{N}\sum_{\delta\in\Delta}\|w_\delta\|^2 =
\frac{|\Delta|}N\ ,
$$
which yields the approximate expression:
\begin{equation}
\Ex{\cD_\Psi(x)} \approx C + \frac{|\Lambda|}N \log_2(\sigma^2)
+ \left(1-\frac{|\Lambda|}N\right)\,
\log_2\left(\tsigma^2\frac{|\Delta|}N\right)\ .
\end{equation}
Therefore, if the ``$\Psi$-component'' of the signal is sparse enough,
i.e. if $|\Lambda|/N$ is sufficiently small (compared with 1),
$\Ex{\cD_\Psi(x)}$ may be expected to behave as 
$\log_2\left(\tsigma^2\frac{|\Delta|}N\right)$, which suggests to
use
\begin{equation}
\hat N_{\psi}(x) = 2^{\cD_{\Psi}(x)}
\end{equation}
as an estimate (up to a multiplicative constant)
for the ``size'' of the $W$ component of the signal.
Notice that this expression coincides with~(\ref{fo:N}),
\item[{\em b.}] The difference between the lower and
upper bounds depends on two parameters: the sparsity
$|\Lambda|/N$ of the $\Psi$-component, and the relative
redundancy parameters $\epsilon(\Delta)$. The latter actually describe
the intrinsic differences between the two considered bases.
When the bases are significantly different, the relative redundancy
may be expected to be small (notice that in any case, it is smaller than 1),
\item[{\em c.}] The relative redundancy parameters $\epsilon$
and $\tepsilon$ which pop up in our model
differs from the one which is generally considered in the literature,
namely the {\em coherence} of the dictionary $W\cup\Psi$
(see e.g.~\cite{Donoho01uncertainty,Elad02generalized,Gribonval03sparse})
$$
\mu[W\cup\Psi]
= \sup_{\stackrel{b,b'\in W\cup\Psi}{b\ne b'}}|\langle b,b'\rangle|\ ,
$$
and the Babel function (see~\cite{Tropp02greed,Gribonval03sparse}.)
The latter are intrinsic to the dictionary, while the Parseval weights
and corresponding $\epsilon$ and $\tepsilon$
provide a finer information, as they also account for the
signal models, via their dependence in the significance maps $\Lambda$ and
$\Delta$,
\item[{\em d.}] Precise estimates for $\epsilon$ and $\tepsilon$
are fairly difficult to obtain\footnote{Our numerical results using
wavelet and MDCT bases suggest that these numbers are generally of
the order of 1/4: any waveform from a given basis always finds
a waveform from the other basis which ``looks like it''.}.
What would actually be needed is a tractable model for the
significance maps $\Delta$ and $\Lambda$, in the spirit of the
structured models described in the two previous sections (for which
we couldn't obtain simple estimates.)
Returning to the wavelet and MDCT case, it is quite natural to
expect that models implementing time persistence in $\Delta$
and scale persistence in $\Lambda$ would yield smaller values for the
relative redundancies than models featuring uniformly distributed
significance maps.
\end{itemize}
A more detailed analysis of this method (including a discussion
of noise robustness issues) is presented in~\cite{Molla03determining}.
\subsection{Numerical simulations}
The above discussion suggest to use the logarithmic dimensions
in order to get estimates for the relative sizes of the
tonal and transient layers in audio signals. We shall use the
following estimated proportions
\begin{equation}
\hat I_{ton} = \frac{\hat N_{\psi}}{\hat N_{\psi}+\hat N_{w}}\ ;\qquad
\hat I_{tr} = \frac{\hat N_{w}}{\hat N_{\psi}+\hat N_{w}}\ ,
\end{equation}

\begin{figure}
\centerline{
\epsfig{file=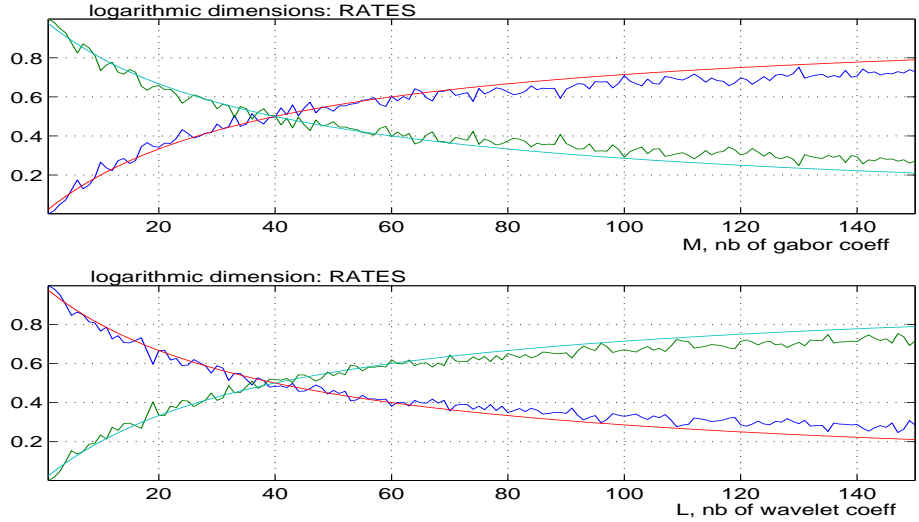, height=7cm,width=12cm}
%\epsfig{file=ThDim25.eps, height=5cm,width=6cm}
%\hskip5mm
%\epsfig{file=ThDim50.eps, height=5cm,width=6cm}
}
\caption{Simulations of tonality and transientness indices, as functions
of $|\Delta|$ or $|\Lambda|$ (time frames of 1024 samples):
theoretical curves and simulation (averaged over 10 realizations);
top plot: $|\Lambda|=40$, varying $|\Delta|$, $I_{tr}$ (decreasing curve)
and $I_{ton}$ (increasing curve);
top plot: $|\Delta|=40$, varying $|\Lambda|$, $I_{ton}$ (decreasing curve)
and $I_{tr}$ (increasing curve.)}
\label{fi:indices.simul}
\end{figure}

In order to validate this approach, we computed these quantities
on simulated signals of the form~(\ref{fo:sig.model}), as functions
of $|\Delta|$ (resp. $|\Lambda|$) for fixed values of $|\Lambda|$
(resp. $|\Delta|$.) The result of such simulations is displayed
in {\sc Figure}~\ref{fi:indices.simul}, which show $\hat I_{ton}$
and $\hat I_{tr}$ as functions of $|\Delta|$, together with the
theoretical curves defined in~(\ref{fo:indices}), averaged over
20 realizations. As may be seen, the results are fairly satisfactory,
which indicates that such indicator may be used for estimating the
percentage of bit rate to be allowed to the different
components, prior to the hybrid coding itself.

\begin{figure}
\centerline{
\epsfig{file=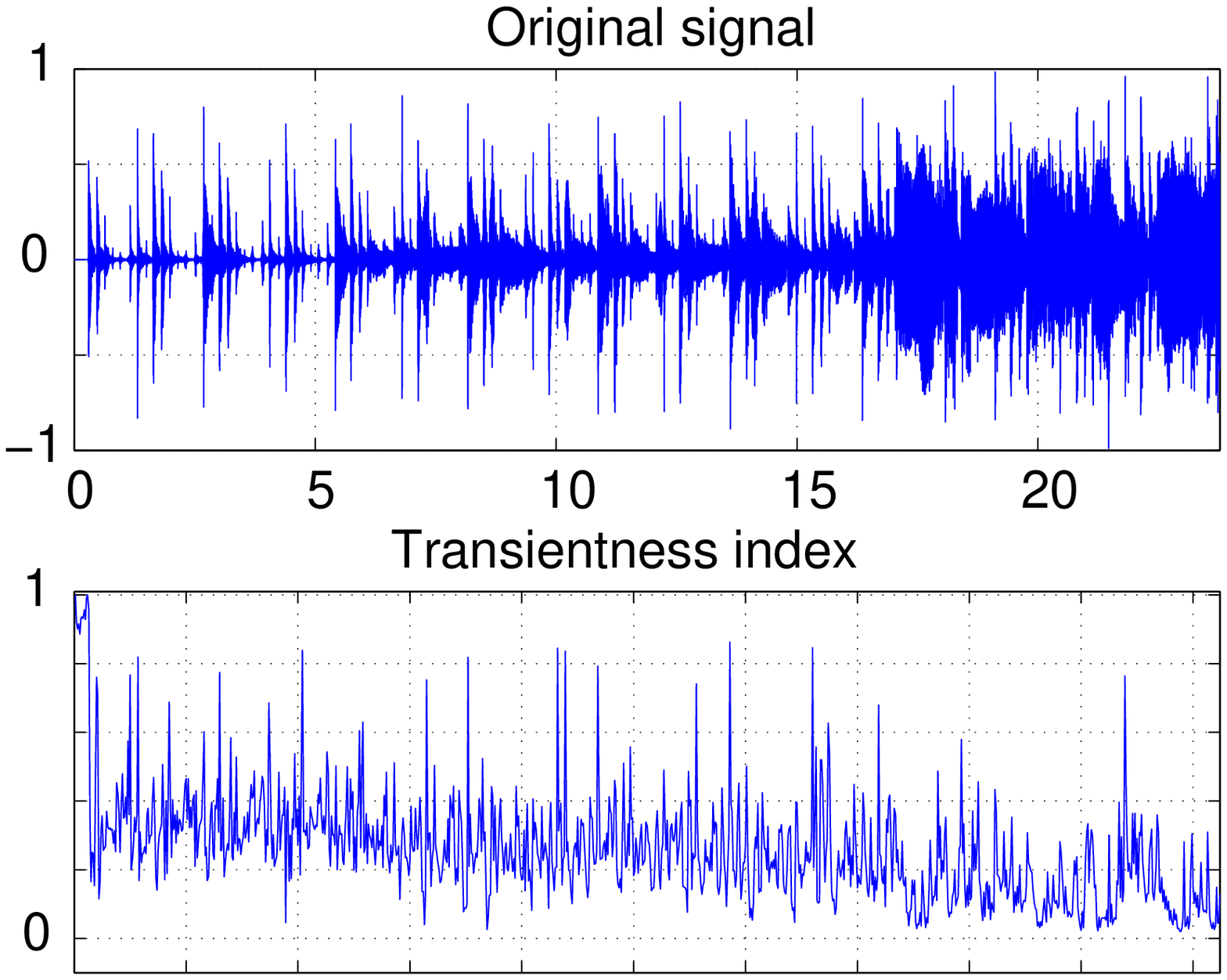, height=8cm,width=5.8cm}\quad
\epsfig{file=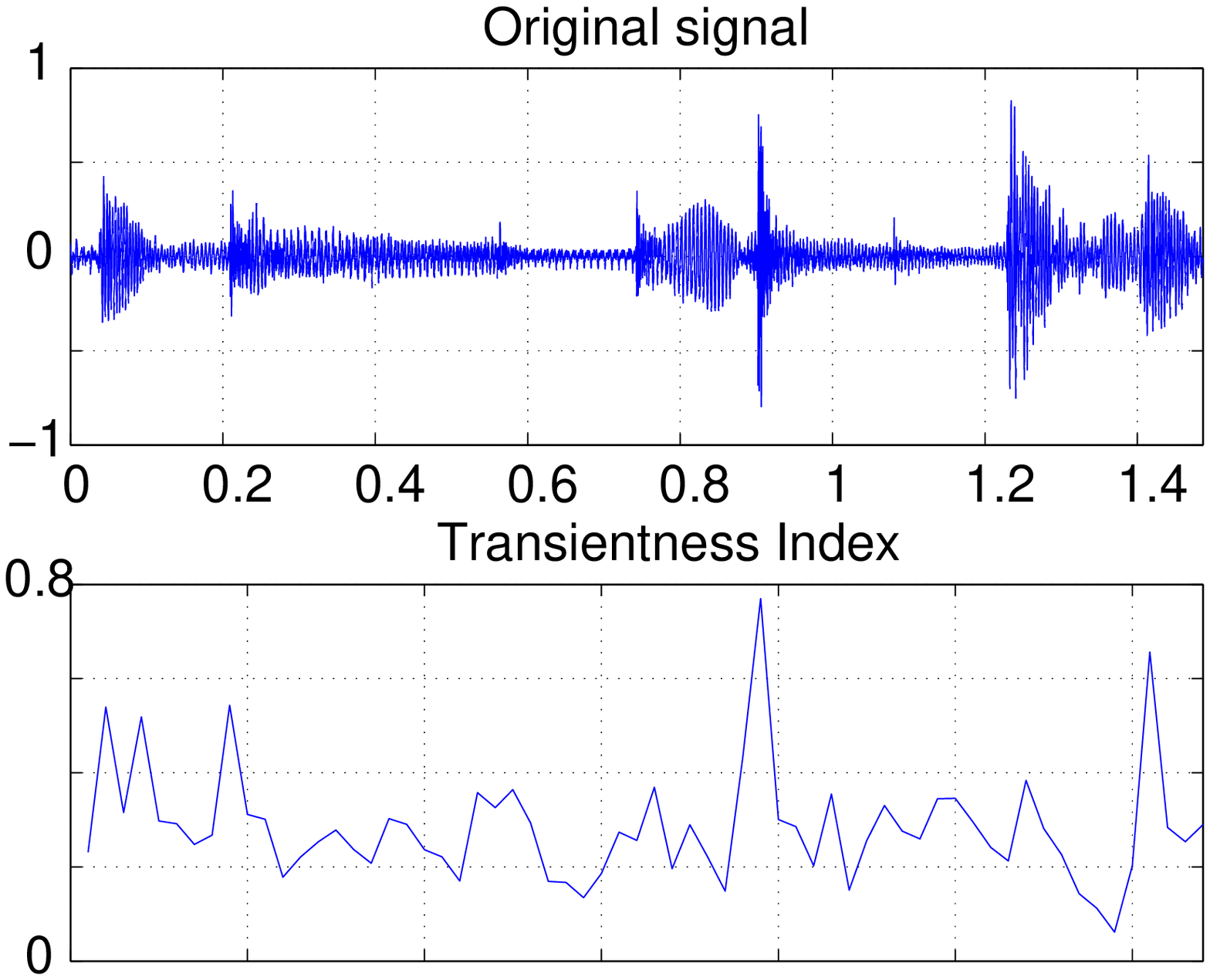, height=8cm,width=5.8cm}
}
\caption{Tonal vs Transient balance for a real audio signal
(a musical signal.) Left: long signal (about 23 seconds long):
top plot: original signal; bottom plot: transientness index.
Right: shorter (1.5 seconds long) segment, same legend.
}
\label{fi:indices.real}
\end{figure}

%\begin{figure}
%\centerline{
%\epsfig{file=sMamaIndex.eps, height=7cm,width=12cm}
%}
%\caption{Transientness index for a shorter (1.5 seconds long) segment of
%the {\em mamavatu} signal.}
%\label{fi:indices.real.2}
%\end{figure}

An example on real audio signal is displayed in
{\sc Figure}~\ref{fi:indices.real},
which represents the transientness index (from which the tonality index
is easily deduced) for a segment (about 23 seconds)
of audio signal (the {\em mamavatu} signal\footnote{
\label{foot:web}available at the web site\\
{\tt http://www.cmi.univ-mrs.fr/$\tilde{ }$torresan/papers/Markov}}, which
will be used again as illustration in the next section.) A shorter
segment of 1.5 seconds (located in the middle of the large segment)
is analyzed similarly in the right hand plots of
{\sc Figure}~\ref{fi:indices.real}
As may be seen, the transientness index (lower curves) exhibits
significant local maxima in the neighborhood of the various
``attacks'' of the signal (see the left hand plots of
{\sc Figure}~\ref{fi:indices.real}.)
Notice also on the right hand plots of {\sc Figure}~\ref{fi:indices.real}
that the transientness
index exhibits an overall decay in the rightmost part of the plot. This
is mainly due to the fact that a significant tonal component shows up in that
part of the signal (see {\sc Figure}~\ref{fi:hyb.exp} in the next section),
which reduces the {\em proportion} of transients
(we recall that the transientness index really measures the proportion, 
and not the {\em quantity} of transient signal present.)

\begin{remark}\rm
It is worth noticing that the indices $\hat I_{ton}$ and
$\hat I_{tr}$ perform satisfactorily as long as the
two expansions in~(\ref{fo:sig.model})
are sparse enough. Otherwise, deviations from the ``ideal''
behavior have to be expected, as may be seen in
the right hand side of the plots in {\sc Figure}~\ref{fi:indices.simul}.
\end{remark}

\begin{remark}\rm
Also, $\hat I_{ton}$ and $\hat I_{tr}$ provide estimates for the
sizes of significance maps only when the variances $\sigma^2$
and $\tsigma^2$ are of comparable magnitude. When this is not the case,
it is easily seen that they rather provide estimates on the
relative energies of the two layers, for example
$\hat I_{tr} = |\Lambda|\sigma^2 /(|\Lambda|\sigma^2 +|\Delta|\tsigma^2)$.
The behavior of the indices in noisy situations (i.e. with small,
additive white noise) may be studied as well, and yields similar
conclusions, as long as the noise's energy is small
enough~\cite{Molla03determining}.
\end{remark}

\section{Application to audio coding}
\label{se:audio}
The ideas developed above are currently being implemented
within a prototype hybrid audio coder, extending the ideas
already described in~\cite{Daudet02hybrid}. While the idea of
hybrid coding of audio signals is not new, our approach
is the first one that implements hybrid transform coding
without prior (time) segmentation of the signal.
A detailed account of the coding system will be given
in a forthcoming publication, together with systematic
performance evaluations. However, we find it interesting
to sketch the main features here, as they provide a thorough
applications of the probabilistic models we just described.

\subsection{Description of a prototype coder}

\begin{figure}[here]
\centerline{
\epsfig{file=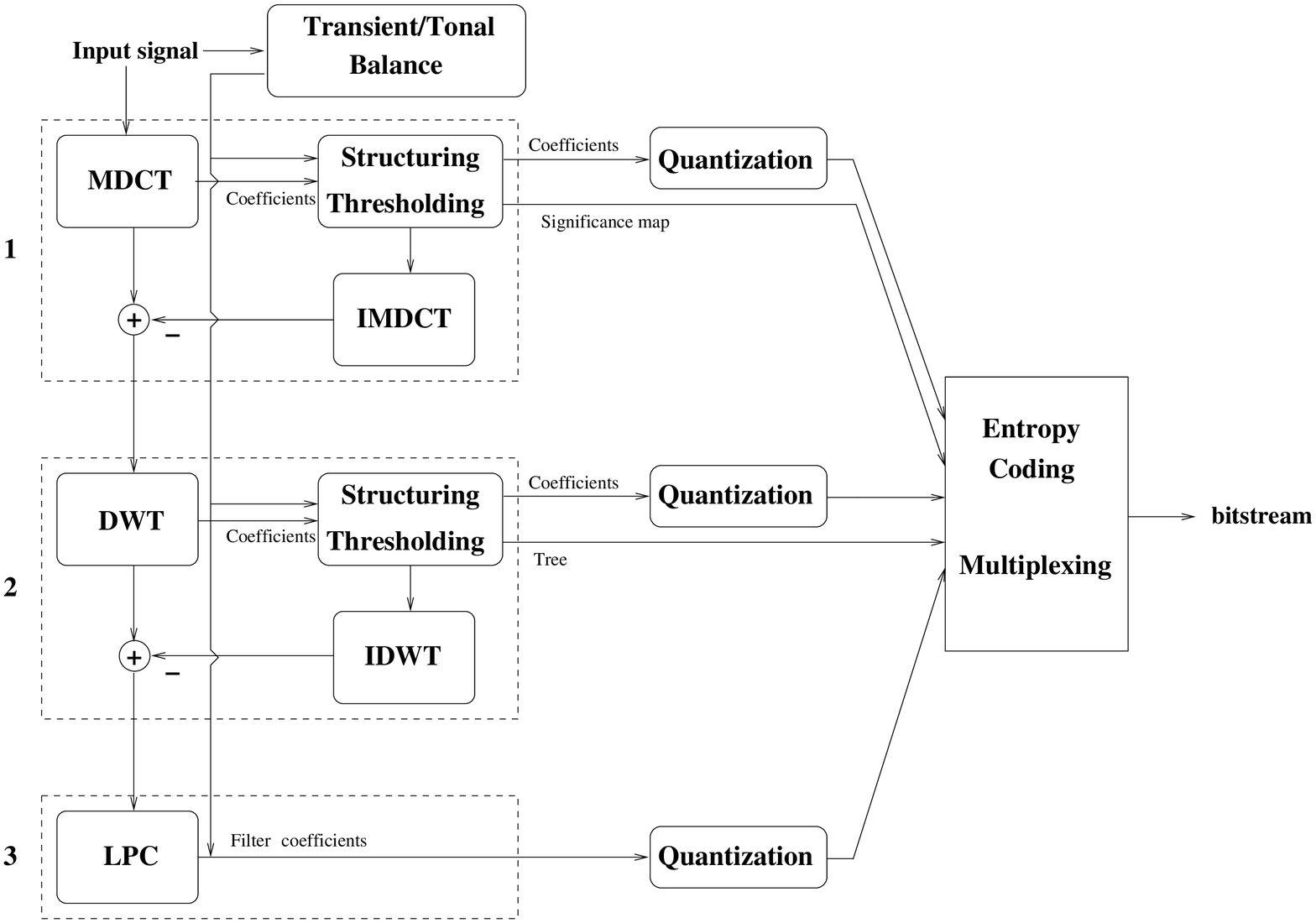, height=12cm,width=12cm}
}
\caption{Block diagram of the hybrid audio coding scheme}
\label{fi:coder}
\end{figure}

The block diagram of the encoder is displayed in
{\sc Figure}~\ref{fi:coder} (the corresponding decoder
simply amounts to invert MDCT and DWT\footnote{DWT: Discrete Wavelet
Transform.} transforms from
the encoded coefficients).
The first step of the algorithm is a pre-estimation of
the relative sizes of the tonal and transient layers,
according to the discussion of section~\ref{se:balance}.
Hence, any given bit budget may be allocated a priori
to the different layers of the signal.

The second step is the estimation of the (structured)
tonal layer, according to section~\ref{se:tonal}.
The parameters of the hidden Markov models are estimated
and updated on large time frames, and the hidden states
are estimated by thresholding of a posteriori probability.
This yields estimated tonal and non-tonal layers
$$
x_{ton} = \sum_{\delta\in\Delta} \langle x,w_\delta\rangle w_\delta\ ;
\quad x_{nton} = x -x_{ton}\ .
$$
The tonal layer is then quantized and encoded using standard
techniques (either uniform quantization, or Lloyd-Max quantization,
for gaussian sources, followed by entropy coding),
while the non tonal layer is transmitted to
the transient layer estimator. Since the parameters of the model
(i.e. the persistence probabilities) provide explicitly
the probabilities of lengths of ``tonal structures'',
the corresponding Huffmann code is readily obtained, and used
for encoding the significance map.

The third step is the estimation of the transient layer
from the non-tonal component. Again, transform coding
is computed within time frames of about 23 milliseconds.
The parameters of the hidden Markov model are estimated,
and updated on larger time frames. Hidden states (i.e.
the significance map) are estimated within each (small)
time frame by thresholding of a posteriori state probability.
Once the transient layer $x_{tr}$ has been estimated,
it is substracted from the signal to yield the residual;
in parallel, the coefficients are quantized and entropy coded.
The tree structure of the transient significance map make it possible
to derive an efficient way of encoding it (see~\cite{Molla03hidden}.)
$$
x_{tr} =
\sum_{\lambda\in\Lambda} \langle x_{nton},\psi_\lambda\rangle\psi_\lambda \ ;
\quad x_{res} = x_{nton} -x_{tr}\ .
$$
The residual is finally modeled as a (locally) stationary
random process, and currently encoded as such using fairly classical
LPC procedures (even though this might not be the optimal solution for
very low bit rate; this subject is currently under study.)

Notice that while the encoding procedure is quite complex (involving
fairly sophisticated estimation algorithms), the decoding is
extremely simple. The tonal and transient layers are reconstructed
on the basis of their significance maps and corresponding encoded
coefficients. The residual is re-generated using LPC technique.

\subsection{Numerical illustrations}
An example of hybrid (or multilayered) signal expansion obtained using
the technique described in this paper is shown
in {\sc Figure}~\ref{fi:hyb.exp} (see {\sc Figure}~\ref{fi:indices.real} for
the corresponding transientness index.) In that example  6\% of coefficients
were retained (no coefficient quantization was done, so this essentially
represents only the ``functional'' part of the compression.) 

\begin{figure}
\centerline{
\epsfig{file=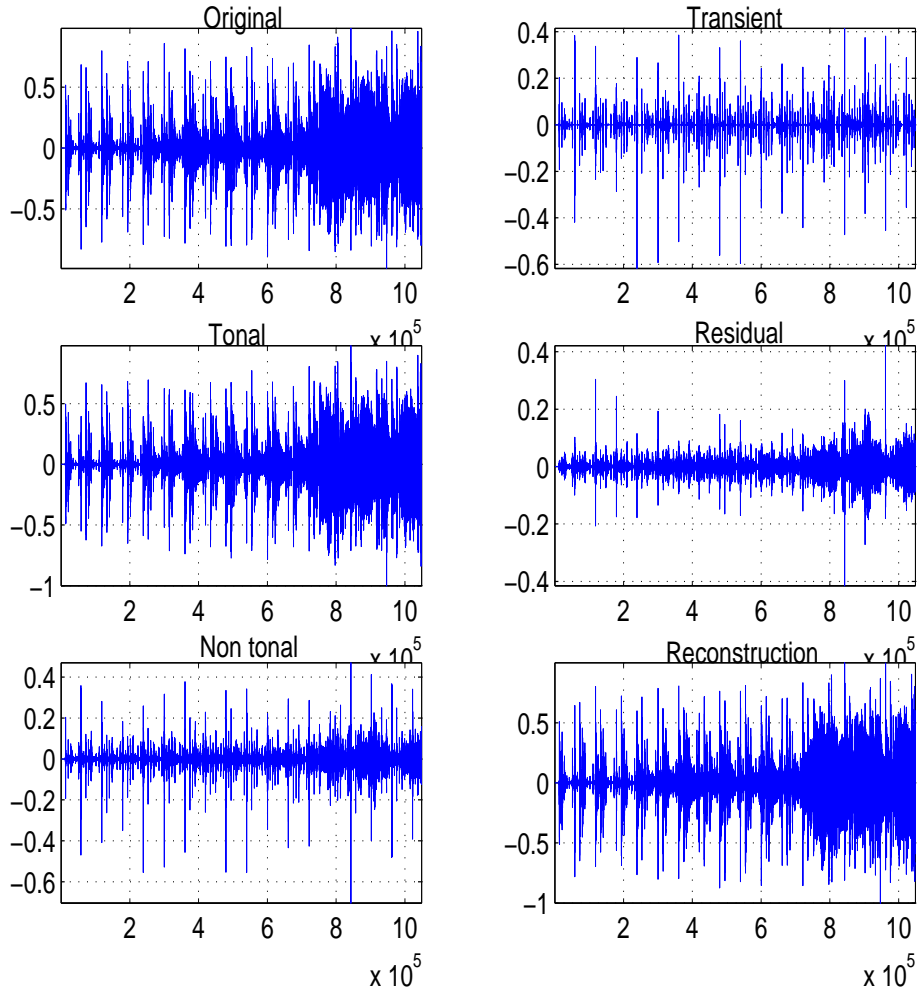, height=14cm,width=12cm}
}
\caption{Compressed hybrid expansion of a piece of musics
(mamavatu, about 23 seconds long.) From top to bottom,
and from left to right: original signal, tonal layer,
nontonal signal, transient layer, residual layer, and reconstruction from
the three layers.}
\label{fi:hyb.exp}
\end{figure}

To demonstrate the ability of the proposed procedure to yield
good signal approximations, we display in {\sc Figure}~\ref{fi:SNR}
a comparison of Signal to Noise (SNR) curves for various encoding
techniques, namely classical MDCT and DWT transform coding,
as well as hybrid coding as proposed in~\cite{Daudet02hybrid}
and the approach described in this paper.
As may be seen, the
performances of the different approaches are essentially
comparable, the standard DWT and MDCT transform coding techniques
(in the left curve) being better by a few dB.
The N-term hybrid method appears to be slightly better than the Markov one.
This is not really surprizing, as
the SNR is computed from the $L^2$ distortion, and the introduction of
structures in the approximation cannot improve the $L^2$ distortion
in comparison with simple coefficient thresholding.
Nevertheless, we notice that this effect is a very weak one.

For the same reason,
introducing a normalization on MDCT coefficients prior to the
selection (right-hand part of {\sc Figure}~\ref{fi:SNR}) strongly penalizes
the N-term MDCT method in terms of $L^2$ distortion. Interestingly
enough, this does not seem to be the case for the two hybrid methods.

\begin{figure}

\centerline{
\epsfig{file=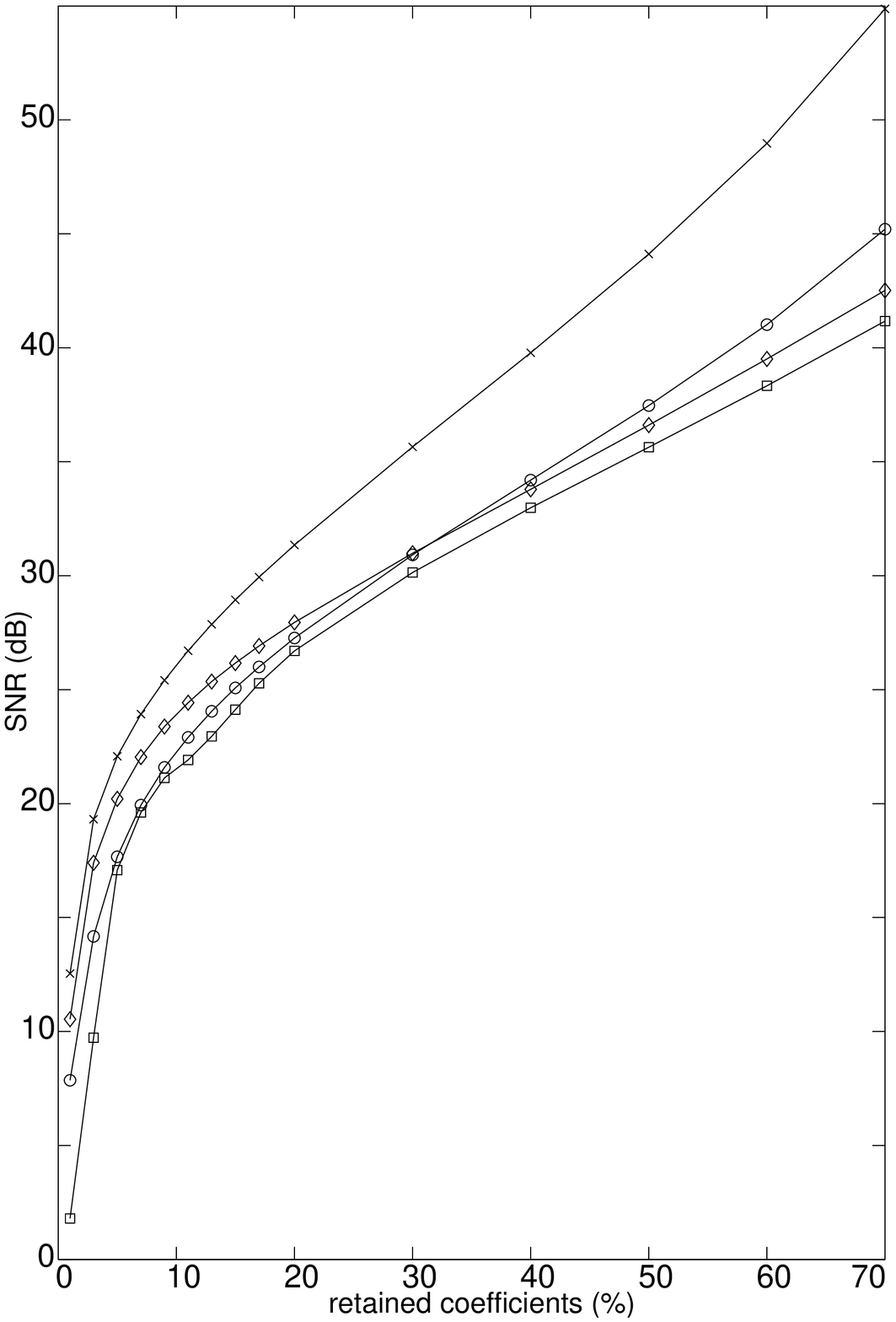, height=8cm,width=8cm}
\epsfig{file=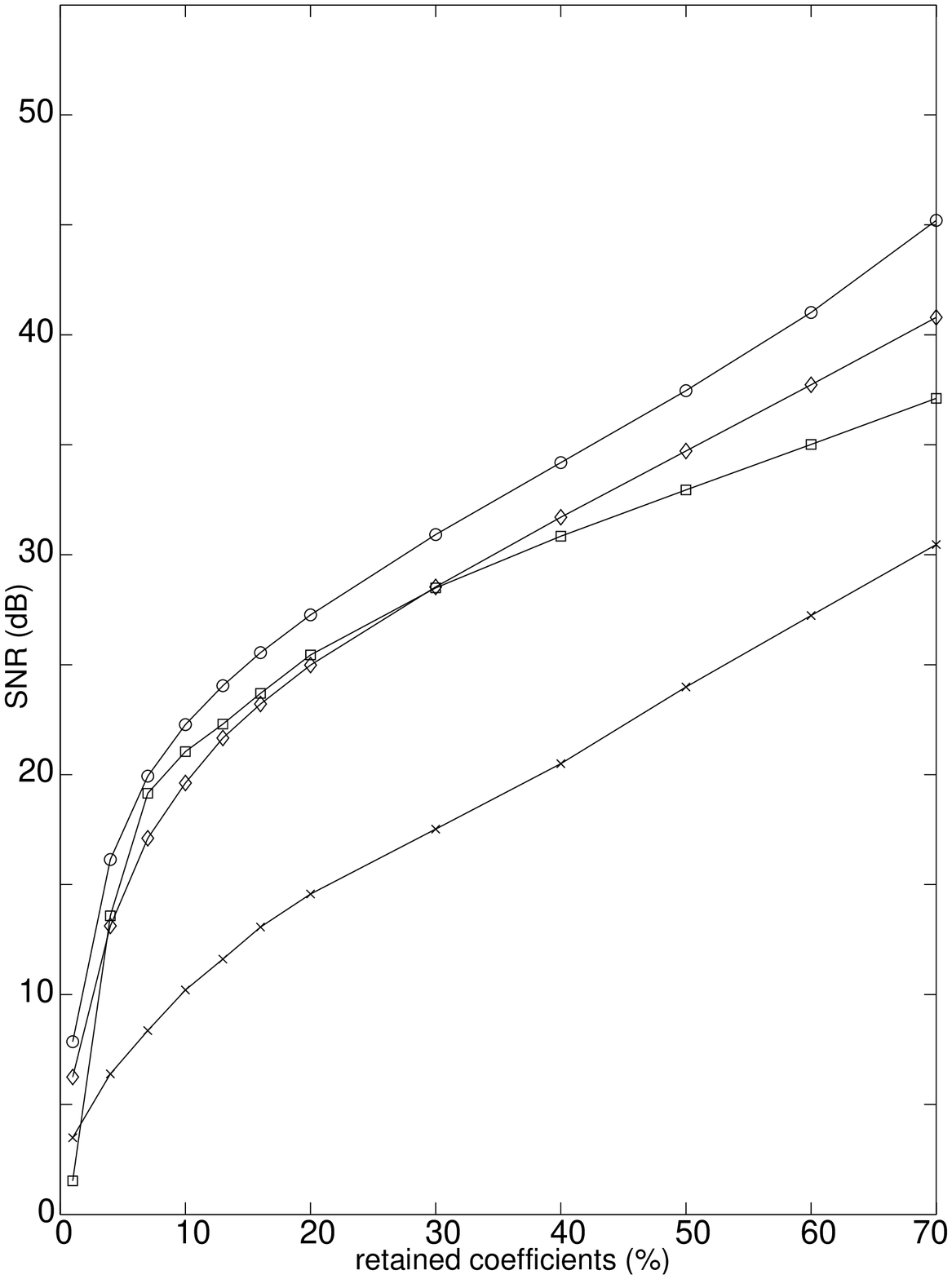, height=8cm,width=8cm}
}
\caption{Comparison of ``functional'' Signal to Noise Ratios for various
coding techniques, computed on the same test signal.
Two configurations are represented: in the left-hand plot, coefficients were
selected without prior normalization in the N-term MDCT (x), DWT ($\circ$)
and hybrid ($\diamond$) methods, and $\alpha$ was equivalently set to 0
in the Markov method ($\square$, see Remark~\ref{rem:freq.decay}).
In the right-hand plot, MDCT coefficients were normalized by a linear
function of the frequency in the N-term MDCT and hybrid methods (no
change in the wavelet method), and $\alpha$ was set to 1 in the
Markov method.}
\label{fi:SNR}
\end{figure}

\medskip
One main outcome of the proposed method is that it provides
a decomposition of the signal into several layers (two, or three if
the residual is considered), i.e. one has to encode several signals
instead of one. Nevertheless, {\sc Figure}~\ref{fi:SNR} (right)
shows that in terms of SNR, this approach remains comparable with
a standard DWT approximation scheme. In addition, the significance maps
are encoded much more efficiently in the Markov approach
(see Proposition~\ref{prop:ton.SMrate} and Corollary~\ref{th:tree_encoding}).

\begin{remark}
\rm
It is worth pointing out an important aspect of audio signal coding.
It is well known that the $L^2$ distortion is {\em not} an adequate
distortion measure for audio signals: it does not take into account
the variation of the hearing threshold as a function of frequency
(which could be done by an appropriate weighting of the $L^2$ norm
in the frequency domain), nor the nonlinear masking effects, which are
extremely important from the perceptual point of view (see for
example~\cite{Painter00perceptual} for a review). The natural distortion
measure that introduces itself naturally in the present scheme
is the likelihood, which takes into account the structures
proposed in the model. For that reason, we believe that such a
distortion measure is more ``natural'' from a perceptual point of view,
though we do not have a clear evidence yet.
\end{remark}

To illustrate the latter remark, the interested reader may find on
the companion web site of this paper (see footnote~\ref{foot:web})
further illustrations of the behaviour of the proposed coding technique
on real audio signals, as well as corresponding sound examples.

\section{Conclusions and perspectives}
\label{se:concl}
This work was based on the belief that efficient signal modeling
cannot be based solely on considerations on individual coefficients
on a well chosen basis (even though one generally tries to
use ``bases that decorrelate'') and may be seen as an
attempt to systematically exploit ``structures'', or ``persistence
properties'' in the coefficient domain. In this respect,
the main contributions of this article are the new hybrid model
we propose, and the a priori rate estimations which may be
deduced from it, thanks to the relative simplicity of the
model (First order Markov chains, and Gaussian distributions.)

We have specially emphasized in Section~\ref{se:audio} the application
to audio coding and compression.
More details on the current implementation of the codec will be
published elsewhere~\cite{Daudet03HSAC} together with a more
complete analysis of quantization issues, and more
detailed numerical results.
Further developments involve designing coefficient quantization procedures
specifically adapted to the tonal and transient layers, as well as
the implementation of adapted masking methods (frequency masking and
time masking.)

Finally, we would also like to point out that compression is
far from being the only application of such models, and
that coding is not limited to compression applications. An efficient
coding scheme, such as the one we propose, should also prove useful
for various applications such as automatic music transcription
(exploiting the tonal layer) and onset detection (exploiting the
transient layer). The multilayered signal representation could also
simplify other audio processing tasks such as time-stretching,
or various other signal modification problems, which may thus be
performed directly in the coefficient domain. Among the
other potential applications of such techniques, let us also mention
blind source separation, which we plan to investigate in the future.

%\section{conclusions}

\bigskip

%\newpage

{\bf Acknowledgements. }
This work was supported in part by the European Union's
Human Potential Programme, under contract HPRN-CT-2002-00285 (HASSIP.)
%S. Molla is supported by the French Ministry of Education
We also acknowledge
support from the AMADEUS Austrian-French exchange programme, which
allowed us to visit the NuHAG group in Vienna, where we had
stimulating exchanges of ideas.
We also wish to thank L. Daudet, F. Jaillet, Ph. Guillemain and
R. Kronland-Martinet for many stimulating discussions.

%\newpage

\bibliographystyle{plain}
\bibliography{Hassip,Markov}

\end{document}